\documentclass[11pt]{article}
\pdfoutput=1
\usepackage{microtype}
\usepackage{lineno}
\usepackage{lmodern}
\usepackage{amsmath,amssymb,amsthm}
\usepackage{xspace}

\usepackage{tikz}
\newcommand{\yibai}{} 
\usepackage{pgfplots}
\usetikzlibrary{calc}
\pgfplotsset{compat=1.13}
\usepgfplotslibrary{fillbetween}

\usepackage{todonotes}

\usepackage{mathtools}

\usepackage[style=alphabetic,maxbibnames=100,maxcitenames=10,backend=bibtex]{biblatex}
\renewbibmacro{in:}{}

\usepackage[colorlinks = true,
			citecolor = blue,
			linkcolor = blue]{hyperref}
\newcommand{\dd}{\mathop{}\!\mathrm{d}}
\let\del\partial

\newcommand{\xaddspace}[0]{\mathchoice{\hspace{-0.8em}}
                {\hspace{-0.6em}}
                {\hspace{-0.4em}}
                {\hspace{-0.3em}}
                 }

\newcommand{\xint}[1]{\int\foreach \i in {2,...,#1}{\xaddspace\int}}
\newcommand{\myiint}[0]{\xint{2}}
\renewcommand{\iint}{\myiint}

\newcommand{\RR}{\mathring R}
\newcommand{\RRR}{\mathring{\bar R}}
\newcommand{\ootimes}{\mathbin{\mathring{\otimes}}}

\let\div\relax
\DeclareMathOperator{\div}{div}
\DeclareMathOperator{\tr}{tr}
\newcommand{\vin}[0]{v^{\textup{in}}}
\newcommand{\vell}[0]{v_{\ell_q}}
\newcommand{\pell}[0]{P_{\ell_q}}
\newcommand{\Rell}{\RR_{\ell_q}}
\newcommand{\pex}[0]{\textsf{\textit{P}}}
\newcommand{\vex}[0]{\textsf{\textit{v}}}
\newcommand{\vv}[0]{\bar v}

\newcommand{\pp}[0]{\bar P}
\newcommand{\wpq}[1][]{w^{\textup{p}}_{#1}}
\newcommand{\wcq}[1][]{w^{\textup{c}}_{#1}}
\newcommand{\wdq}[1][]{w^{\textup{d}}_{#1}}

\newcommand{\Rlinear}{R_{q+1}^\textup{linear}}
\newcommand{\Rosc}{R_{q+1}^\textup{osc}}
\newcommand{\TT}[0]{\mathsf{T}}
\newcommand{\Id}[0]{\textup{Id}}
\newcommand{\glue}[0]{\mathsf g}

\newcommand{\sobolev}[2]{{W^{#1,#2}}}
\newcommand{\mylinear}[1]{\mathbf{#1}}
\newcommand{\lv}[0]{\mylinear{v}}
\newcommand{\lw}[0]{\mylinear{w}}
\newcommand{\lP}[0]{\mylinear{P}}
\newcommand{\lf}[0]{\mylinear{f}}
\newcommand{\jmax}[0]{j_{\max}}
\newcommand{\imax}[0]{i_{\max}}
\newcommand{\mykappa}[0]{{\kappa}}
\newcommand{\mymu}[0]{{\mu}}
\newcommand{\myK}[0]{{K}}

\def\dashint{\,\ThisStyle{\ensurestackMath{%
  \stackinset{c}{.2\LMpt}{c}{.5\LMpt}{\SavedStyle-}{\SavedStyle\phantom{\int}}}%
  \setbox0=\hbox{$\SavedStyle\int\,$}\kern-\wd0}\int}
\def\ddashint{\,\ThisStyle{\ensurestackMath{%
  \stackinset{c}{.2\LMpt}{c}{.5\LMpt+.2\LMex}{\SavedStyle-}{%
    \stackinset{c}{.2\LMpt}{c}{.5\LMpt-.2\LMex}{\SavedStyle-}{%
      \SavedStyle\phantom{\int}}}}\setbox0=\hbox{$\SavedStyle\int\,$}\kern-\wd0}\int}

\usepackage{scalerel,stackengine}
\stackMath
\newcommand\widecheck[1]{%
\savestack{\tmpbox}{\stretchto{%
  \scaleto{%
    \scalerel*[\widthof{\ensuremath{#1}}]{\kern-.6pt\bigwedge\kern-.6pt}%
    {\rule[-\textheight/2]{1ex}{\textheight}}
  }{\textheight}%
}{0.5ex}}%
\stackon[1pt]{#1}{\scalebox{-1}{\tmpbox}}%
}

\newcommand{\coloneq}{\mathrel{\mathop:}=}
\newcommand{\eqcolon}{=\mathrel{\mathop:}}

\newcommand{\ee}{\textup{e}}
\DeclareMathOperator{\curl}{curl}
\newcommand{\kk}{\overline  k}
\newcommand{\kkk}{k'}
\newcommand{\bk}{{(k)}}
\newcommand{\bj}{{(j)}}
\newcommand{\TR}[1]{\widetilde{R_{#1}}}

\newcommand{\varEPS}{0.2}
\newcommand{\varTAU}{0.2}
\newcommand{\varNi}{2}
\newcommand{\varN}{4}
\newcommand{\varYOFF}{-8*\varTAU }
\newcommand{\varETASTART}{\varNi*\varTAU}

\newcommand{\colorR}{blue!50!white}
\newcommand{\colorEta}{blue!20!white}
\newcommand{\colorEtaMinusOne}{red}

\newcommand{\crvi}[1]{
    \varTAU/3 - (1-2*\varEPS)*\varTAU/3 
    }
\newcommand{\crvii}[1]{
    2*\varTAU/3 + (1-2*\varEPS)*\varTAU/3 
    }

\newcommand{\crviii}[1]{
    \varTAU/3 - (1-\varEPS)*\varTAU/3 
    + 2*\varEPS *\varTAU *sin(deg(2 * pi * #1))/3
    }
\newcommand{\crviv}[1]{
    2*\varTAU/3 + (1-\varEPS)*\varTAU/3 
    + 2*\varEPS *\varTAU *sin(deg(2 * pi * #1))/3
    }
    
\newcommand{\etaminusone}[1]{
    (1 - (#1)*(#1) ) * (1 - (#1)*(#1))
}

\newtheorem{thm}{Theorem}[section]
\newtheorem{cor}[thm]{Corollary}

\newtheorem{lem}[thm]{Lemma}
\newtheorem{prop}[thm]{Proposition}
\theoremstyle{definition}

\theoremstyle{remark}
\newtheorem{rem}[thm]{Remark}
\numberwithin{equation}{section}
\newcommand{\Holder}{H\"older\xspace}

\newcommand{\CZ}{Calder\'on--Zygmund\xspace} 

\DeclareMathOperator{\supp}{supp}
\bibliography{convexbib}
\newcommand{\betain}[0]{{\beta^{\textup{in}}}}%

\begin{document}
\title{Construction of solutions to the 3D Euler equations with initial data in $H^\beta$ for $\beta>0$}
\author{Calvin Khor and Changxing Miao}
\maketitle
\begin{abstract}
In this paper, we use the method of convex integration to construct infinitely many distributional solutions in $H^{\beta}$ for $0<\beta\ll1$ to the initial value problem for  the three-dimensional incompressible Euler equations. We show that if the initial data has any small fractional derivative in $L^2$, then we can construct solutions with some regularity, so that the corresponding $L^2$ energy is continuous in time. This is distinct from the $L^2$ existence result of  
\textit{E. Wiedemann}, Ann. Inst. Henri Poincar\'e, Anal. Non Lin\'eaire 28, No. 5, 727--730 (2011; Zbl 1228.35172), where the energy is discontinuous at $0$.
\\[0.4em]
\end{abstract}

\section{Introduction}

\label{s:intro} 
We study the Cauchy problem for the 3D incompressible Euler equations,
\begin{equation}
\left\{ \begin{alignedat}{-1}
\del_t v+(v\cdot\nabla) v  +\nabla P   &=  0,
 \\
  \nabla \cdot v  &= 0,
\end{alignedat}\right.  \label{e:euler}
\end{equation}
 which describes the motion of inviscid fluids. $v(x,t)\in\mathbb R^3$ is the fluid velocity, and $P(x,t)\in\mathbb R$ is the pressure field that maintains the incompressibility condition. Although this classical partial differential equation has been studied for many years, a global existence result for low regularity data only appeared  relatively recently in 2011: 
by appealing to a result of De Lellis and Sz\'ekelyhidi \cite{zbMATH05675020}, Wiedemann  \cite{zbMATH05965634} constructed   global weak solutions to the 3D Euler equations for an arbitrary periodic (weakly) divergence-free initial data in $L^2$. Furthermore, these solutions resulted from the method of convex integration and are therefore highly non-unique. The solutions constructed are bounded in time with values in $L^2_\sigma$, and continuous in time with values in $L^2_\sigma$ with the weak topology, i.e. they belong to $L^\infty L^2_\sigma \cap C^0 L^2_{\sigma,\textup{weak}}$. In addition, the energy of these solutions  are discontinuous at the initial time.

Convex integration is also the main tool of our paper, and has roots in the work of Nash \cite{zbMATH03093832} which showed the existence of rather arbitrary $C^1$ isometric embeddings of a Riemannian manifold. The connection of this geometric problem to fluid mechanics was discovered by De Lellis and Sz\'ekelyhidi \cite{zbMATH05710190}. The method starts with some sort of `subsolution' and iteratively cancels the difference from a true solution by adding carefully selected oscillatory terms.  We refer the reader to the survey articles \cite{zbMATH07207865,zbMATH07301372,zbMATH07310889} for an overview of the technique.

The particular method used in  \cite{zbMATH05965634} has since been refined, leading to the resolution \cite{zbMATH06976275} of Onsager's conjecture: distributional $\beta$-\Holder Euler flows must conserve energy if $\beta>1/3$, and counterexamples exist for $\beta<1/3$. Later, Buckmaster and Vicol \cite{zbMATH07003146} found extra freedom in the convex integration procedure to take advantage of `intermittency', which is a saturation of the Bernstein inequalities in the convex integration building blocks. Our aim is to see what these new developments bring to the Cauchy problem for the 3D Euler equation in $H^\beta$ spaces. 

For any given initial data in $H^\betain$, we construct  $H^\beta$ weak solutions (with $0<\beta \ll \betain$) of the 3D Euler equations on the spatially periodic domain $\mathbb T^3 = (\mathbb R/ \mathbb Z)^3$:

\begin{thm}\label{t:main}
Let $T\in[0,\infty]$, $\betain>0$ and $\vin \in H^{\betain} (\mathbb T^3)$. Then for $\beta< \min(\frac\betain{216}, \frac1{66})$, there exists infinitely many smooth functions $e(t): [0,T]\to [0,\infty)$ with $e|_{[1,T]} > \|\vin\|_{L^2}^2$, and a corresponding  solution $v\in C^0([0,T]; H^\beta) $ to \eqref{e:euler} on $[0,T]$ with initial data $\vin$, such that for all $t\in [1,T]$,
\[\int_{\mathbb T^3} |v(x,t)|^2 \dd x = e(t).\]
\end{thm}
We give a more precise version of Theorem \ref{t:main} below (Theorem \ref{t:precise}) after defining some parameters.

Our constant $\frac1{66}\approx 2^{-6.04}$ is rather small, though it is a little larger than the current best constant for Navier--Stokes  $\approx 2^{-16}$ in \cite{zbMATH07003146}. (We expect a mild improvement is possible for Navier--Stokes but leave this to future work.) We obviously cannot expect non-uniqueness for $\beta>5/2$, since this would imply they solve the equation pointwise with a Lipschitz gradient, and such solutions are unique. It has in fact been shown \cite{zbMATH05291954} that energy is conserved for $\beta>5/6$, which means that our method would need to be drastically altered to go past $5/6$. Such questions are anyway very far out of reach: 
 recently, a long paper of Buckmaster,  Masmoudi,  Novack and  Vicol \cite{https://doi.org/10.48550/arxiv.2101.09278} (see also \cite{https://doi.org/10.48550/arxiv.2203.13115}) showed that there exist energy non-conservative solutions in $C^0 H^{1/2-\epsilon}$, but their intricate `pipe dodging' seems to be at odds with Isett's gluing technique. We crucially use gluing to allow arbitrary initial data, so we are currently unable to take advantage of the methods developed there. 

As mentioned earlier, the strong ill-posedness for $L^2$ data was shown by Weidemann \cite{zbMATH05965634}. Scheffer showed many years prior \cite{zbMATH00427755} that the trivial solution was distributionally non-unique in $L^2$. \cite{Rosa2021DimensionOT} showed (in particular) that all smooth data can give rise to non-unique solutions. For data in \Holder spaces, the work of Daneri, Sz\'ekelyhidi and Runa \cite{zbMATH07370998} showed that a dense set of $C^{1/3-\epsilon}$ data give rise to multiple solutions. Their solutions are also dissipative, i.e. energy non-increasing. From weak-strong uniqueness  \cite{zbMATH06999781}, we know that such initial data are not smooth. 

In contrast, the solutions that we construct  must \emph{increase the energy}. The authors of the present paper with Ye have recently \cite{https://doi.org/10.48550/arxiv.2204.03344} shown non-uniqueness for \emph{all} data in $C^{1/3-\epsilon}$ by allowing the energy to rise. The present paper can be thought of as a continuation of this work, although the calculations are different and the proof idea differs in some key areas, which we discuss further in Subsection \ref{ss:comparison}. The energy increase is also seen in \cite{Rosa2021DimensionOT} and \cite{zbMATH05965634}. 

Finally, we  also mention some of the ill-posedness results for the Euler equations of the norm-inflation type, i.e. ill-posedness in a certain functional framework: see \cite{zbMATH06468723,zbMATH06540654,zbMATH06919596}. 
There is also the possibility of finite-time blow-up, which has been rigourously shown to occur \cite{zbMATH07441733, https://doi.org/10.48550/arxiv.1910.14071}.
\section{Outline of the convex integration scheme}
\label{s:outline}
We consider the `Euler--Reynolds' system $(v_q,P_q,\RR_q)$ with initial data,
\begin{equation}
\left\{ \begin{alignedat}{-1}
\del_t v_q+\div (v_q\otimes v_q)  +\nabla P_q   &=  \div \RR_q,
 \\
  \nabla \cdot v_q &= 0,
  \\ v_q |_{t=0}&=\vin* \psi_{\ell_{q-1}},
\end{alignedat}\right.  \label{e:subsol-euler}
\end{equation}
where $\ell_{q-1}$ is defined below in \eqref{e:ell},
$\psi_{\ell_{q-1}}(x):=\ell_{q-1}^{-3} \psi(\ell_{q-1}^{-1}x)$ is a  standard mollifier at scale $\ell_{q-1}$, $a \otimes b\coloneq ab^\TT  = (a_i b_j)_{i,j=1}^3$, and the divergence of a matrix $M=(M_{ij})_{i,j=1}^3$ is the vector $\div M$ defined as
\begin{align*}
(\div M)_i \coloneq  \sum_{j=1}^3 \partial_j M_{ij}.
\end{align*}
In particular, since $v_q$ is divergence-free, $\div (v_q\otimes v_q) = (v_q \cdot\nabla )v_q$.
 The tensor $\RR_q$ is required to be a symmetric, trace-free $3\times3$ matrix, i.e.
\begin{gather}
     \RR_q = \RR_q^\TT ,  \qquad \tr\RR_q = \sum_{i=1}^3 (\RR_q)_{ii} =  0. \label{e:RR-sym-tracefree-cond}
\end{gather}
We enforce the mean zero conditions on $v_q$ and $P_q$:
 \begin{align}
     \int_{\mathbb T^3} v_q  \dd x= 0, \qquad \int_{\mathbb T^3} P_q \dd x = 0. \label{e:mean-zero}
 \end{align}
 These conditions are natural from the Gallilean symmetries of the equation.
\subsection{Parameters and the main iterative proposition}
\label{ss:params}
For given integers $a,b$ satisfying $a>1$, $b>1$, and $0<\betain<1$, $\beta<\min(\frac1{10},\betain)$,
and a given $\vin\in H^{\betain}$ with
\begin{align}
        \|\vin\|_{H^\betain}&\le 1,\label{e:vin-CK}
\end{align} we define
\begin{align}
    \lambda_q \coloneq   a^{b^q} ,\qquad  \delta_q &\coloneq \lambda_q^{-2\beta}. \notag%
\end{align}
We also define\footnote{This is essentially chosen for \eqref{e:vell-vq} below to hold.} the small positive parameter $\ell_q\le \lambda_q^{-1}$  by
\label{ss:mollification}
\begin{align}
    \ell_q \coloneq \frac{\delta_{q+1}^{1/2}}{\delta_q^{1/2} \lambda_q^{1+3\alpha/2}} = \lambda_q^{-1-3\alpha/2-\beta(b-1)} \label{e:ell}.
\end{align}
Note that \eqref{e:vin-CK} implies that we do not take arbitrary initial data $\vin$, but this is easily fixed in the proof of Theorem \ref{t:main} by a scaling argument.
So without loss of generality, we firstly choose a smooth function $e:[0,\infty]\to (\|\vin\|_{L^2}^2,\infty)$  such that
\begin{align}
t\in[1-\tau_0,T]\implies \delta_{2}\lambda_1^{-\alpha} &\le e(t) -\int_{\mathbb T^3} |\vin|^2 \dd x \le \delta_{2} .
    \label{e:vin-energy-estimate}
\end{align}

The estimates we propagate inductively are:
\begin{align}
    \|v_q\|_{L^2} &\le 1-\delta_q^{1/2} ,
    \label{e:v-q-H0}
    \\
    \|v_q\|_{H^1} &\le M \delta_{q}^{1/2} \lambda_q ,
    \label{e:v-q-H1}
    \\
    \|\RR_q\|_{L^1} &\le \delta_{q+1}\lambda_q^{-3\alpha} ,
    \label{e:RR-q-L1}
    \\
   t\in[ 1-\tau_{q-1},T] \implies \lambda_q^{-\alpha}\delta_{q+1} &\le e(t) -\int_{\mathbb T^3} |v_q(x,t)|^2 \dd x \le \delta_{q+1}.
    \label{e:energy-q-estimate}
\end{align}
The constant $M=361$ is universal and computed in Corollary \ref{c:estimates-for-wpq-wcq}. In particular $M$ does not depend on $a,b,\beta,\alpha$ or $q$. 

\begin{prop}
\label{p:main-prop}Let $\vin$ be a $H^{\betain}(\mathbb T^3)$ function such that $\|\vin\|_{H^\betain}\le 1$, and let $e:[0,T]\to [0,\infty)$ be any smooth function satisfying \eqref{e:energy-q-estimate}.
There exist $\beta_0<\betain$, $b_0=b_0(\beta,\betain )$, $b\ge b_0$  a multiple of six,  $\alpha_0=\alpha_0(\beta,\betain,b)< 1 $, and $a_0=a_0(\beta,\betain,b,\alpha_0)>1$ such that for $0<\alpha < \alpha_0$ and $a>a_0$, the following holds:

Let $(v_q,P_q,\RR_q)$ solve 
\eqref{e:subsol-euler} and satisfy the conditions \eqref{e:RR-sym-tracefree-cond}, \eqref{e:mean-zero}, \eqref{e:v-q-H0}--\eqref{e:energy-q-estimate}.
Then there exist smooth functions $(v_{q+1}, p_{q+1}, \RR_{q+1})$ solving \eqref{e:subsol-euler} and satisfying the conditions \eqref{e:RR-sym-tracefree-cond}, \eqref{e:mean-zero}, \eqref{e:v-q-H0}--\eqref{e:energy-q-estimate}
with $q$ replaced by $q+1$, and such that
\begin{align}
        \|v_{q+1} - v_q\|_{L^2} + \frac1{\lambda_{q+1}}\|v_{q+1}-v_q\|_{H^1} &\le \frac M2 \delta_{q+1}^{1/2}.
        \label{e:velocity-diff-prop2-1-final-est}
\end{align}
\end{prop}
The constants  $b_0$ and $\beta_0$ are explicitly computed in \eqref{e:constraint-intermittency}, \eqref{e:constraint-b}, and \eqref{e:constraint-b-beta-strong}. This shows in particular that we can take $b=6$ and $\beta_0=\min(\frac{\betain}{216},\frac1{66})$.

Our main result is the following theorem, which implies Theorem \ref{t:main}:
\begin{thm}
   \label{t:precise} 
   Let $T\in[0,\infty]$, $\betain>0$ and $\vin \in H^{\betain} (\mathbb T^3)$. Choose  $\beta<\beta_0,b>b_0,\alpha<\alpha_0$ and $a>a_0$ as in Proposition \ref{p:main-prop}. Also define $\Gamma$ as in \eqref{e:gamma} below.
        Then for any smooth function  $e(t): [0,T]\to [0,\infty)$  solving the inequality 
\begin{align}
 \Gamma^2 \delta_2 \lambda_1^{-\alpha} \le e(t) - \int_{\mathbb T^3}|\vin|^2 \dd x \le \Gamma^2 \delta_2,\quad t\in [1,T], \label{e:rescaled-energy-ineq}
\end{align}
there exists a corresponding weak solution $v\in C^0([0,T]; H^\beta) $ to \eqref{e:euler} on $[0,T]$ with initial data $\vin$, such that for all $t\in [1,T]$,
\[\int_{\mathbb T^3} |v(x,t)|^2 \dd x = e(t).\]
\end{thm}
\begin{rem} The upper and lower bounds of \eqref{e:rescaled-energy-ineq} explode as e.g. $a\to\infty$. A better result can be obtained for special initial data e.g. if $\vin$ is itself a stationary Euler flow. It would be interesting (as also in \cite{https://doi.org/10.48550/arxiv.2204.03344}) to improve this to bounds that are independent of certain parameters, or to find ways to constrain the energy for all times $t\ge 0$.
\end{rem}

\subsection{Proof of Theorem \ref{t:precise} using Proposition \ref{p:main-prop}}

Choose parameters $b,\beta,\betain,\alpha,a$ and initial data $\vin$ with $\|\vin \|_{H^\betain}\le 1$ as in Proposition \ref{p:main-prop}. Define 
\begin{gather*}
    v_1\coloneq \vin*\psi_{\ell_{0}},
    \qquad  p_1\coloneq |v_1|^2-\int_{\mathbb T^3} |v_1|^2\dd x,
    \\ \RR_1\coloneq v_1\ootimes v_1  \coloneq  v_1\otimes v_1 -|v_1|^2 I_{3\times 3}.
\end{gather*}
 Assume first that \eqref{e:v-q-H0}, \eqref{e:v-q-H1} and \eqref{e:RR-q-L1} hold for $q=1$. We  choose $e$ so that \eqref{e:vin-energy-estimate} holds and is otherwise arbitrary.
 Applying Proposition \ref{p:main-prop}, we obtain a sequence of functions $v_q$ that solve \eqref{e:subsol-euler} and satisfy the conditions \eqref{e:RR-sym-tracefree-cond}, \eqref{e:mean-zero}, \eqref{e:v-q-H0}--\eqref{e:energy-q-estimate}. In particular, they converge in $C^0_t L^2_x$ as $\sum_{q=1}^\infty  \| v_{q+1} - v_q\|_{L^2} < \infty$, and in fact in $C^0 H^{\tilde \beta}$ for all $\tilde \beta < \beta$, as
 \[ \sum_{q=1}^\infty \|v_{q+1} - v_q\|_{H^{\tilde \beta}} \lesssim  \sum_{q=1}^\infty \|v_{q+1} - v_q\|_{L^2}^{1-\tilde \beta} \|v_{q+1} - v_q\|_{H^1 }^{\tilde \beta}  =\sum_{q=1}^\infty \delta_{q+1} \lambda_{q+1}^{\tilde \beta} < \infty.  \]
 
This proves Theorem \ref{t:precise}, under the assumption that $\|\vin\|_{L^2_x} \le 1$, and that \eqref{e:v-q-H0},\eqref{e:v-q-H1} and \eqref{e:RR-q-L1} hold. For a general initial data in $H^\betain$, we define the rescaling
\begin{gather}
 \Gamma = \max\left(1,\frac{\|\vin\|_{L^2}}{1-\delta_1^{1/2}}, \frac{\|\vin\ootimes \vin\|_{L^1}^{1/2} }{\delta_{2}^{1/2} \lambda_1^{-3\alpha/2} }, \frac{(1+\ell_0^{-1})\|\vin\|_{L^2}}{M\delta_1^{1/2}\lambda_1}\right), \label{e:gamma} \\ 
\widetilde{\vin} \coloneq \vin/\Gamma, \qquad \widetilde{v_1} \coloneq \widetilde{\vin} * \psi_{\ell_0}. \label{e:vq,1-tilde}
\end{gather}
By construction, $\widetilde{v_1}$ solves \eqref{e:v-q-H0}, \eqref{e:v-q-H1}, and \eqref{e:RR-q-L1}. So  we can apply the earlier construction together with any function $\widetilde e$ satisfying \eqref{e:vin-energy-estimate} with $\Gamma T$ in place of $T$. This produces a weak solution $\widetilde v$ to \eqref{e:euler} such that $\widetilde v(x,0)= \widetilde \vin$, and that $\|\widetilde v(t)\|_{L^2}^2 = \widetilde e(t) > \|\widetilde{ \vin}\|_{L^2}^2$ for $t\ge 1$. To finish, we simply undo the scaling transform by setting
\[ v(x,t) \coloneq  \Gamma v(x,\Gamma t), \quad \text{and} \quad e(t) \coloneq \Gamma^2 \tilde e (\Gamma t). \]
Then $v$ is also a solution to the Euler equation. For times $t\ge 1/\Gamma$ (and hence for times $t\ge1$),  its energy solves \eqref{e:rescaled-energy-ineq}, so in particular $\|v(t)\|_{H^2}^2 = e(t)>\|\vin\|_{L^2}^2$. \hfill \qedsymbol\\

The remainder of the paper is devoted to the proof of Proposition \ref{p:main-prop}. 

\subsection{Sketch for  Proposition \ref{p:main-prop} and organisation of paper}
The proof follows the same scheme of other convex integration papers that has been optimised over the years (e.g. see \cite{zbMATH07038033}): starting from the triplet $(v_q, P_q, \RR_q)$, in Section \ref{s:construct}, we suitably prepare them for the main convex integration step by a mollification and then Isett's gluing procedure \cite{zbMATH06976275}. Then in Section \ref{s:wq+1-defn}, we use some intermittent Mikado flows as in \cite{zbMATH07301372} (although we only use spatial intermittency; see the definition of $W_\bk$), with some modified cutoffs from our earlier work \cite{https://doi.org/10.48550/arxiv.2204.03344} to control the energy at positive times. In Section \ref{s:new-stress}, we define and estimate the resulting new stress error $\RR_{q+1}$.  In Section \ref{s:energy-iteration}, we demonstrate the required energy control \eqref{e:energy-q-estimate} with $q+1$ in place of $q$.



\subsection{Comparison with earlier work {\cite{https://doi.org/10.48550/arxiv.2204.03344}}}
\label{ss:comparison}
As can be seen in the proof of Theorem \ref{t:precise}, we need to propagate estimates for $v_q$ in $L^2$. It follows that the natural space for the  error $\RR_q$ which arises from the quadratic nonlinearity is $L^1$. This is the cause for the differences between this paper and \cite{https://doi.org/10.48550/arxiv.2204.03344}, where \Holder estimates are used for both $v_q$ and $\RR_q$. 


One such difference is the size of certain estimates. The gluing parameter $\tau_q$ defined in \eqref{e:tau_q-and-t_i} performs the same function as in Isett's  work, but is much smaller as a function of $\lambda_q$. The reason for this difference is that the well-posedness requires at least $5/2$ derivatives in the $L^2$ scale, but only 1 derivative in the \Holder scale. Another difference which is somewhat minor is that we resort to $L^\infty$ estimates when estimating $L^p$ norms of compositions, to make use of Fa\`a di Bruno estimates. Improving this estimate would lead to a minor improvement in the constants that we do not pursue in this paper. 

These considerations lead to the following heuristic which should aid readers of this paper:  for most building blocks which are not compositions of functions, a gradient costs $\ell_q^{-1}$, and a time derivative costs $\tau_q^{-1}$. For other functions (e.g. see Proposition \ref{p:chi}), the Fa\`a di Bruno estimate makes gradients cost $\ell_q^{-4}$, and a time derivative costs $\tau_q^{-1} \ell_q^{-3}$.

The other main difference is how we achieve the initial data in our construction. In \cite{https://doi.org/10.48550/arxiv.2204.03344}, in order to achieve the $C^{1/3-\epsilon}$ regularity, we incrementally improved the initial data at the gluing step. The gluing step is not as well-behaved in the $L^2$ case, and in any case what we can prove is far from what should be the optimal regularity. Hence, in this paper, we simply adjust the initial data by adding the stationary flow
\[ \wdq \coloneq  \vin * \psi_{\ell_q} - \vin * \psi_{\ell_{q-1}}. \]
 For our choices of constants, and crucially for data in $H^\betain$ with $\betain>0$,  this is a very small term (see Corollary \ref{c:estimates-for-wpq-wcq}) and does not significantly disturb the other estimates.

\label{ss:sketch-pf-main-prop}

\section{Initial preparations: mollification and gluing}
In this section, we suitably modify $v_q$ so that we can perform convex integration with Mikado flows in Section \ref{s:wq+1-defn}. 

In the sequel, we will use the notation $A\lesssim B$ to mean that $A\le C B$ for a constant that may change from line to line, but does not depend on the parameters $a,b,$ or $q$. We will also write $A\ll B$ to mean that $A\lesssim CB$ where $C$ is some negative power of $\lambda_q$. For such inequalities, by making $a$ sufficiently large, we can achieve the inequality $A\le B$ with constant 1.

\label{s:construct}

\subsection{Mollification}
Observe that,  if
$
\beta(b-1) < \frac1{K_0}, 
$ and  $\alpha\ll1$, then  $ \lambda_q^{-1-1/K_0} \le  \ell_q$. If we choose $b=6$ and $\beta<1/66$ then this is true with $K_0=13$, which is much better than the estimates we use below.
The functions $\vell$ and $\Rell$ are defined as
\begin{align}
    \vell &\coloneq v_q * \psi_{\ell_q},
 &   \Rell &\coloneq \RR_q * \psi_{\ell_q}  + \Big(\vell \ootimes \vell - (v_q \ootimes v_q) * \psi_{\ell_q}  \Big) , \label{e:vell-Rell}
\end{align}
where $\ell_q$ is defined in \eqref{e:ell}, and for  $a,b\in\mathbb R^3$, we have written  $a\ootimes b \coloneq a\otimes b - (a\cdot b)\Id$ for a trace-free version of $a\otimes b$.  $\vell$ is also divergence-free, and $\Rell$ is also symmetric and trace-free. They satisfy the equations
\begin{equation}
\left\{ \begin{alignedat}{-1}
\del_t \vell +\div (\vell\otimes \vell)  +\nabla \pell   &=  \div \Rell ,
\\
  \nabla \cdot \vell &= 0,
\end{alignedat}\right.  \label{e:mollified-euler}
\end{equation}
where $\pell \coloneq P_q *\psi_{\ell_q} -|v_q|^2 + |\vell|^2,$ and we have used the identity $\div(f\Id) =\nabla f$ for scalar fields $f$. The following proposition follows from standard mollification estimates, and the following  Constantin--E--Titi identity from \cite{zbMATH06082181}: For any multiindex $\alpha$, and any two constants $f_0,g_0$, $f_0=f_0*\psi_\epsilon$, $g_0=g_0*\psi_\epsilon$  and $(f_0g_0)*\psi_\epsilon =(f_0*\psi_\epsilon )(g_0*\psi_\epsilon)$; it follows that
\begin{align}
     \MoveEqLeft \del^\alpha\Big ( (f* \psi_{\epsilon}) (g* \psi_{\epsilon})- (f\otimes g)*\psi_{\epsilon}\Big) (x) \notag
      \\&= \Big((f-f_0)(g-g_0)\Big) * \del^\alpha \psi_\epsilon  (x) \label{e:C-E-T}
      \\&\quad  - \sum_{0\le \delta \le \alpha} \binom \alpha\delta (f-f_0)*\del^\delta\psi_\epsilon (x) (g-g_0)*\del^{\alpha-\delta} \psi_\epsilon (x) .  \notag
\end{align}
In particular, this is true for each fixed $x$, with  $f_0=f(x)$, $g_0=g(x)$.
\begin{prop}[Estimates for mollified functions]\label{p:estimates-for-mollified} For all $N\ge 0$,
\begin{align}
\|\vell\|_{L^2} &\le 1-\delta_q^{1/2} \label{e:vell-L^2} , 
\\ 
 \quad \|\vell\|_{H^{N+1}} &\le M \delta_{q}^{1 / 2} \lambda_{q}\ell_q^{-N},  \label{e:vell-HN+1}
\\
\|\vell\|_{C^{N}}&\lesssim \delta_q^{1/2}\lambda_q \ell_q^{-N-1/2-\alpha}, \label{e:vell-Cn}
\\
\|\Rell\|_{W^{N,1}} &\lesssim \delta_{q+1} \lambda_q^{-3\alpha} \ell_q^{-N} , \label{e:Rell}
\\
\|\vell-v_{q}\|_{L^2} &\lesssim \delta_q^{1/2}\lambda_q\ell_q =\delta_{q+1}^{1 / 2}  \lambda_{q}^{-3\alpha/2}, \label{e:vell-vq}
\\
\Big|\int_{\mathbb{T}^{3}} | v_{q}|^{2}-|\vell|^{2} \dd x \Big| &\lesssim \delta_{q+1}  \lambda_q^{-3\alpha} .\label{e:energy-vell}
\end{align}
\end{prop}
\begin{proof}
   \eqref{e:vell-L^2} is because mollifiers do not increase $L^p$ norms. Since $|\vell(x) - v_q(x)| \le  \int_{|y|\le\ell_q} |y| \psi_{\ell_q}(y) \int_0^1 |\nabla v_q(x-ty)| \dd t \dd y$, the $H^1$ estimate \eqref{e:v-q-H1} gives \eqref{e:vell-vq}. Putting $N$ derivatives on the mollifier gives \eqref{e:vell-HN+1}. \eqref{e:vell-Cn} comes from  \eqref{e:vell-HN+1} and Sobolev embedding $H^{3/2+\alpha}\subset C^0$. 
    \eqref{e:energy-vell} follows from \eqref{e:C-E-T} with $\alpha=0$ because $\int_{\mathbb T^3} |v_q|^2 \dd x = \int_{\mathbb T^3} |v_q|^2 * \psi_{\ell_q} \dd x$. Finally, \eqref{e:Rell} follows by using \eqref{e:C-E-T} with $|\alpha|\le N$ to estimate the commutator term in \eqref{e:vell-Rell}.
\end{proof}
\subsection{Classical exact  flows}
\label{ss:exact}
We define  $\tau_q\in(0,\ell_q^{\alpha})$ and the  sequence of initial times $t_i$ ($i\in\mathbb Z_{\ge0}$) by
 \begin{align}
    \tau_q \coloneq \frac{\ell_q^{\frac{3}{2}+2\alpha}}{\delta_q^{1/2} \lambda_q}, \qquad t_i\coloneq i\tau_q,   \label{e:tau_q-and-t_i}
\end{align}
and  define  $(\vex_i,\pex_i)$ for $i\ge 0$ to be the unique (exact) solutions to the Euler equations, defined on some interval around $t_i$, with initial data
 $\vex_i|_{t=t_i} = \vell(t_i)$. That is,
  \begin{align}
 \left\{ \begin{alignedat}{-1}
\del_t \vex_i +(\vex_i\cdot\nabla) \vex_i  +\nabla \pex_i   &=  0,
\\
  \nabla \cdot \vex_i  &= 0,
  \\
  \vex_i|_{t=t_i} &= \vell(\cdot,t_i)
\end{alignedat}\right. & \label{e:exact-euler}
\end{align}

\begin{prop}[Energy Estimates {\cite[Section 3.2.3]{zbMATH01644218}}]
\label{p:exact-euler}
Let $\tilde s> \frac{5}{2}$. The unique $H^{\tilde s}$ solution $V$ to \eqref{e:euler} with initial data $V_0\in H^{s}(\mathbb T^3)$ is defined at least for $t\in [-\tilde T,\tilde T]$, where $\tilde T= c\|\nabla V_0\|^{-1}_{L^\infty_x }$ for some universal $c>0$, and satisfies the following  estimates for $N\in [\frac{5}{2},\tilde N]$ and all $t\in[-\tilde T,\tilde T]$,
\[ \|V(\cdot,t)\|_{L^2_x}=\|V_0\|_{L^2_x},\qquad  \|V(\cdot,t)\|_{H^{N}_x} \lesssim  \| V_0 \|_{H^{N}_x} . \]
\end{prop}
Note that the definition of $\tau_q$ was chosen\footnote{note $\|\vell\|_{H^1}\tau_q\lesssim \|\vell\|_{H^{\frac{5}{2}+\alpha}}\tau_q \lesssim \ell_q^{\alpha}$,
 so $\|\vell\|_{H^1} \le \tau_q^{-1}\ell_q^{\alpha}.$} so that (with $a\gg 1$ to absorb the constant $c$) $\vex_i$ is well-defined on $[t_{i-1},t_{i+1}]$, since $ H^{3/2+\alpha}\subset L^\infty$. In particular, we have
\begin{cor}\label{p:exact-euler-cor-HN-estimates}
    $(\vex_i,\pex_i)$ are defined on $t_i+[-\tau_q,\tau_q]$, and for all $N\ge 0$,
     \[ \|\vex_i\|_{L^2} \le  \|\vell\|_{L^2}\le 1-\delta_q^{1/2},\quad  \|\vex_i\|_{H^{N+1}}\lesssim \|\vell\|_{H^{N+1}} \lesssim \delta_q^{1/2} \lambda_q\ell_q^{-N}. \]
\end{cor}
\begin{proof}
This follows from the standard $H^s$ energy estimate for any $s\ge 0$, as derived in e.g. \cite{zbMATH03915851},
$ \frac{\dd}{\dd t}\|\vex_i\|_{H^s}^2 \lesssim \|\nabla \vex_i\|_{L^\infty} \|\vex_i\|_{H^s}^2 $,
and the fact that $\int_{t_i-\tau_q}^{t_i+\tau_q}\|\nabla \vex_i \|_{L^\infty}\dd t \ll   1 $ by the choice of $\tau_q$ and \eqref{e:vell-Cn}.
\end{proof}
\begin{prop}[Stability]
\label{p:stability}
For each $i\ge 0$, $|t-t_i|\lesssim  \tau_q$, $p\in(1,\infty)$, and $N\ge 0$, 
\begin{align}
    \|\vex_i -\vell \|_{W^{N,p}} &\lesssim  \tau_q \delta_{q+1}\lambda_q^{-3\alpha}\ell_q^{-N-1-3/p'} \label{e:stability-v-WNp}.
    \\
        \|\del_t (\vex_i -\vell) \|_{W^{N,p}} &\lesssim    \delta_{q+1}\lambda_q^{-3\alpha}\ell_q^{-N-1-3/p'} \label{e:stability-v-Dt-WNp}.
\end{align}\end{prop}
\begin{proof} 
This follows from  standard $L^p$ \emph{a priori} estimates which we now sketch. 
 The difference $\lv:=\vex_i - \vell$ satisfies the equation \eqref{e:abstracted-equation} with
\begin{gather*} \lv = \vex_i - \vell,\qquad  \lw = \vell,\qquad  \lP = \pex_i - \pell, \\\lf = -\div \Rell - \lv\cdot\nabla \vex_i \eqcolon f_{0,1}+f_{0,2}.\end{gather*} 
Since $\|f_{0,1}\|_{L^p} \le \|\Rell\|_\sobolev 1p \lesssim \|\Rell\|_{\sobolev {1+3/p'}1}$ and $\|f_{0,2}\|_{L^p} \le \|\nabla \vex_i \|_{L^\infty}\|\lv\|_{L^p_x}$, we obtain from  the definition of $\tau_q$, \eqref{e:abstracted-Lp-estimate}, Corollary \ref{p:exact-euler-cor-HN-estimates} and Gr\"onwall's inequality that
\begin{align*}
     \|\vex_i(\cdot,t) - \vell(\cdot,t) \|_{L^p_x} &\lesssim \underbrace{\ee^{\int_{t_i-\tau_q}^{t_i+\tau_q}\|\nabla v_\ell\|_{L^\infty} \dd s}  }_{\lesssim 1} \int_{t_i-\tau_q}^{t_i+\tau_q}  \|\Rell\|_\sobolev {1+3/p'}1 \dd t' \\
     &\lesssim \tau_q  \|\Rell\|_\sobolev {1+3/p'}1.
\end{align*}

This combined with the estimates \eqref{e:Rell} (and interpolation)  prove  \eqref{e:stability-v-WNp} for $N=0$. 
For the higher derivatives one proceeds via induction. For any multiindex $\alpha$ with $|\alpha|=N\ge1$, differentiating the equation for $\vex_i-\vell$ leads to 
\eqref{e:abstracted-equation} with 
\begin{gather}  \lv = D^\alpha (\vex_i - \vell),\qquad  \lw = \vell,\qquad  \lP = D^\alpha(\pex_i - \pell),   \label{e:abstracted-eqn-for-WNp} 
\\
\lf = -D^\alpha \div \Rell - D^\alpha((\vex_i-\vell)\cdot\nabla\vex_i) + (\lw\cdot \nabla \lv  -D^\alpha(\lw\cdot \nabla (\vex_i-\vell)  )\notag  \end{gather} 
Some calculation with Leibniz rule and \eqref{e:abstracted-Lp-estimate} eventually leads to the estimate
\[ \begin{array}{r}\|\vex_i(t) - \vell(t)\|_{W^{N,p}_x} \le \int_{t_i-\tau_q}^{t_i+\tau_q}\bigg( \| \Rell \|_\sobolev {N+1+3/p'}1 + \| \vell\|_{C^N} \|\vex_i - \vell\|_{W^{1,p}_x} \\
 + \|\vex_i - \vell\|_{L^p_x} \| \vex_i\|_{C^{N+1}} + (\|\vex_i\|_{C^1} + \|\vell\|_{C^1})\|\vex_i - \vell\|_{W^{N,p}_x} \bigg)\dd t'. \end{array}\]
From here, the systematic application of estimates and Gronwall's inequality as in the $N=0$ case inductively proves \eqref{e:stability-v-WNp} for all $N\ge1$.

For the second estimate \eqref{e:stability-v-Dt-WNp},  the equation \eqref{e:abstracted-eqn-for-WNp} implies the bound
\[ \|\del_t \lv\|_{L^p} \lesssim \|\lw\|_{L^\infty}\|\nabla \lv\|_{L^p} + \|\nabla \lw\|_{L^\infty}\|\lv\|_{L^p} +\|\lf\|_{L^p}.\]
This can be controlled using \eqref{e:vell-Cn} and \eqref{e:stability-v-WNp}, and the earlier estimates.
\end{proof}

\subsection{Vector potential estimate}
\label{ss:vector-potentials}
By the Helmholtz decomposition for smooth functions on $\mathbb T^3$, a divergence-free field $V$ of zero mean is a curl, i.e. $V=\nabla\times Z$ for an incompressible field  $Z\eqcolon \mathcal BV$  called the vector potential of $V$. The operator $\mathcal B=(-\Delta)^{-1}\curl$ is the `Biot--Savart operator'.

Later, to estimate the glued stress tensor $\RRR_q$ in Proposition \ref{p:estimate-vvq-RRRq}, we will need estimates of $\mathcal R(\vex_i-\vex_{i+1})$, where $\mathcal R$ is the inverse divergence operator recalled in Section \ref{s:inverse-div-op}. By the $L^p$-boundedness of \CZ operators  \cite[Chapter VII]{zbMATH03367521}, these are provided by the following estimates of $z_i-z_{i+1}$:
\begin{prop} \label{p:stability-z}
   Define $z_i \coloneq \mathcal B \vex_i$. For any $i\ge 0$, $|t-t_i|\lesssim  \tau_q$, $p\in(1,\infty)$, and $N\ge 0$,  
\begin{align}\left\|z_{i}-z_{i+1}\right\|_{W^{N,p}} &\lesssim \tau_{q} \delta_{q+1} \lambda_q^{-3\alpha} \ell_q^{-N-3/p'},\label{e:stability-z-WNp}
\\
\left\|\del_t(z_{i}-z_{i+1})\right\|_{W^{N,p}} &\lesssim \delta_{q+1} \lambda_q^{-3\alpha} \ell_q^{-N-3/p'}.\label{e:stability-z-Dt-WNp}
\end{align}
\end{prop}
\begin{proof}To estimate $z_i - z_{i+1}$, it suffices to estimate each $z_i - z_\ell$ where $z_\ell\coloneq \mathcal B \vell $. For $N\ge1$, this directly follows from Proposition \ref{p:stability} and the $L^p$-boundedness of $\nabla\mathcal B$. For $N=0$, substituting $\vex_i-\vell=\nabla\times (z_i-z_{\ell_q})$ in the equation for $\vex_i-\vell$, we find after some  work (as in \cite[Section 3.3]{zbMATH07038033}) the equation \eqref{e:abstracted-equation} with\footnote{We have written $(v\times \nabla)^k$ for the operator $\sum_{k,l,m=1}^3\epsilon_{klm} v^l \del_{x^m}$, where $\epsilon_{klm}$ is the Levi--Civita symbol.}
\begin{gather*} \lv = z_i - z_\ell,\qquad  \lw = \vell,\qquad  \lP = 0, \\\lf = -(-\Delta)^{-1} \bigg[ \nabla\times \div\RR_\ell +\sum_{j,k=1}^3\curl\del_{x^k} \Big(  (\lv \times \nabla)^k \lw^j + (\lv \times \nabla)^j \lw^k \Big)\bigg].\end{gather*} 
\eqref{e:stability-z-WNp} now follows from \eqref{e:abstracted-Lp-estimate}, Gr\"onwall's inequality, and the $L^p$-boundedness of \CZ operators. This estimate also implies \eqref{e:stability-z-Dt-WNp} using the equation, similarly to the proof of \eqref{e:stability-v-Dt-WNp}.
\end{proof}

\subsection{Gluing procedure}
\label{ss:gluing}
Define the intervals $I_i,J_i$ ($i\ge -1)$ by \begin{align}
    I_i &\coloneq \Big[ t_i + \frac{\tau_q}3,\ t_i+\frac{2\tau _q}3\Big] ,& 
     J_i &\coloneq \Big(t_i - \frac{\tau_q}3,\ t_i+\frac{\tau _q}3\Big).
     \label{e:I_i-J_i-defn}
\end{align}
and define 
  $\imax\coloneq\sup\{ i\ge0 :   (J_i\cup I_i)\cap [0,T]\neq \emptyset \}\le \left\lceil \frac T{\tau_q}\right\rceil.$ Their union $\bigcup_{i=0}^{\imax} I_i \cup J_i$ covers $[0,T]$.
Also let $\{ \glue_i\}_{i=0}^{\imax}$  be a partition of unity such that
\begin{align}
     \supp \glue_i &=I_{i-1} \cup J_i \cup I_{i}, & \glue_i |_{J_i} &=1, & \|  \del_t^N  \glue_i\|_{C^0_t} &\lesssim \tau_q^{-N}\ \   (N\ge 0). \label{e:chi_i-properties}
\end{align}
In particular, for $|i-j|\ge2$, $\supp \glue_i \cap \supp \glue_i =\emptyset$. By a slight abuse of notation, we write $ \glue_i \vex_i \coloneq 0$ for $t$ outside the support of $\glue_i$ even when $\vex_i$ is not defined, and similarly for other functions.  Then we define the glued velocity and pressure $\vv_q$ by
\begin{align}
    \vv_q(x,t) &\coloneq \sum_{i=0}^{\imax}  \glue_i(t) \vex_i(x,t).   \label{e:vv_q}
\end{align}
Note that $\bar{v}_q(0,x)=v_{\ell}(0,x)=\vin* \psi_{\ell_{q-1}}* \psi_{\ell_q}$, $\vv_q$ is  divergence-free, and \begin{align*}
    t\in I_i\implies\vv_q(x,t)&=\vex_i(x,t)\glue_i(t) + \vex_{i+1}(x,t)\glue_{i+1}(t),\\
    t\in J_i\implies \vv_q (x,t)&=\vex_i(x,t).
\end{align*}
We therefore define
\begin{align*}
    \RRR_q &\coloneq
     \sum_{i=0}^{\imax} \del_t \glue_i \mathcal R(\vex_i-\vex_{i+1} ) - \glue_i(1-\glue_i)(\vex_i-\vex_{i+1} )\ootimes (\vex_i-\vex_{i+1} ),
    \\
     \pp_q &  \coloneq  \sum_{i=0}^{\imax}\Bigg(  \glue_i(t) \pex_i(x,t) -
      \glue_i(1-\glue_i)\bigg( |\vex_i - \vex_{i+1}|^2  - \int_{\mathbb T^3} |\vex_i - \vex_{i+1}|^2 \dd x\bigg)\Bigg),
\end{align*}
i.e. for $t\in J_i$, $\RRR_q \mathrel{\mathop:}\equiv  0 ,\ \pp_q\mathrel{\mathop:}\equiv \pex_q$, and for $t\in I_i$,
\begin{align}
    \RRR_q &\coloneq
        \del_t \glue_i \mathcal R(\vex_i-\vex_{i+1} ) - \glue_i(1-\glue_i)(\vex_i-\vex_{i+1} )\ootimes (\vex_i-\vex_{i+1} ), \label{e:RRR_q-def} \\
    \pp_q &  \coloneq \glue_i  \pex_q - \glue_i(1-\glue_i)\Big( |\vex_i - \vex_{i+1}|^2  - \int_{\mathbb T^3} |\vex_i - \vex_{i+1}|^2 \dd x\Big), \label{e:def:pP_q}
\end{align}
where we have used the inverse divergence operator $\mathcal R$ as recalled in Section \ref{s:inverse-div-op}. It is easy to check that $\vex_i-\vex_{i+1}$ has mean zero (so that $\RRR_q$ is well-defined), and furthermore, $\pp_q$ has mean zero, while $\RRR_q$ is symmetric and trace-free. The functions $(\vv_q, \pp_q, \RRR_q)$ solve:
\begin{equation}
\left\{ \begin{alignedat}{-1}
\del_t \vv_q+\div (\vv_q\otimes \vv_q)  +\nabla \pp_q   &=  \div \RRR_q,
\\
  \nabla \cdot \vv_q &= 0.
\end{alignedat}\right.  \label{e:subsol-glued-euler}
\end{equation}
\begin{prop}[Estimates for $\vv_q$] \label{p:estimate-vvq-RRRq} For all $N\ge 0$, and all $p\in(1,\infty)$,
\begin{align}
\|\vv_{q}-\vell\|_{L^2} 
& 
\lesssim\tau_q\delta_{q+1}\lambda_q^{-3\alpha}\ell_q^{-5/2} 
= \delta_{q+1}^{1/2}\lambda_q^{-3\alpha/2}\ell_q^{2\alpha} 
\ll  \delta_{q+1}^{1 / 2} \ell_q^\alpha,\label{e:vv_q-stability-L2} 
\\
\|\vv_{q}-\vell\|_\sobolev{N}{p} & \lesssim \tau_{q} \delta_{q+1} \lambda_q^{-3\alpha} \ell_q^{-N-1-3/p'}, \label{e:vv_q-stability-WNp} 
\\
\|\vv_{q}\|_{H^{1+N}} & \lesssim \delta_{q}^{1 / 2} \lambda_{q} \ell_q^{-N},\label{e:vv_q-bound}
\\
\|\RRR_q\|_{W^{N,p}}
    &\lesssim \delta_{q+1} \lambda_q^{-3\alpha} \ell_q^{-N-3/p'}. \label{e:RRR_q-N+alpha-bd}
\\
\|\del_t \RRR_q\|_{W^{N,p}}
    &\lesssim \tau_q^{-1}\delta_{q+1} \lambda_q^{-3\alpha} \ell_q^{-N-3/p'}. \label{e:RRR_q-N+alpha-bd-dt}
    \end{align}
    In addition, we have $p=1,\infty$ estimates with a small loss:
\begin{align}
    \|\RRR_q\|_{W^{N,1}} &\lesssim \delta_{q+1} \lambda_q^{-2\alpha} \ell_q^{-N}\ll \delta_{q+1}\ell_q^{-N+\alpha }, \label{e:RRR_q-WN1}  
\\   
   \|\RRR_q\|_{C^{N}} &\lesssim \delta_{q+1} \lambda_q^{-2\alpha} \ell_q^{-N-3}\ll \delta_{q+1}\ell_q^{-N-3+\alpha }, \label{e:RRR_q-CN}   
\\   
   \|\del_t\RRR_q\|_{C^{N}} &\lesssim \tau_q^{-1} \delta_{q+1} \lambda_q^{-2\alpha} \ell_q^{-N-3}\ll \tau_q^{-1}  \delta_{q+1}\ell_q^{-N-3+\alpha }, \label{e:RRR_q-CN-dt}     
\end{align}

\end{prop}
\begin{proof}
    The  estimates \eqref{e:vv_q-stability-L2}--\eqref{e:RRR_q-N+alpha-bd-dt} follow directly from the $W^{N,p}$ estimates of Propositions \ref{p:stability} and \ref{p:stability-z}. For \eqref{e:RRR_q-WN1}, we use the embedding $W^{N,p} \subset \sobolev {N}{1}$ on the compact space $\mathbb T^3$ for $0<p-1\ll 1$ depending only on $\alpha$. Using $\ell_q \ge \lambda_q^{-3/2}$, it suffices to take $1/p > 1-2\alpha/9$. For \eqref{e:RRR_q-CN} and \eqref{e:RRR_q-CN-dt}, we use this value of $p$ and the Sobolev embedding $W^{N+3,p}\subset C^N$.
\end{proof}
\begin{prop}[Energy of $\vv_q$]
\label{p:energy-of-vv_q}
\begin{align}
   \sup_t \left| \int_{\mathbb T^3} |\vv_q(x, t)|^2-|\vell(x, t)|^2\dd x \right|\lesssim \delta_{q+1}\ell_q^{2\alpha}. \label{e:energy-of-vv_q}
\end{align}
\end{prop}
\begin{proof}
In each $J_i$, $\vv_q=\vex_i$, so that  $\int_{\mathbb T^3} |\vv_q(x,t)|^2\dd x=\int_{\mathbb T^3} |\vex_i(x,t)|^2\dd x$.
For $t\in I_i$, since $\vv_q=\glue_i\vex_i+(1-\glue_i)\vex_{i+1}$, we have
\begin{align*}
    & \int_{\mathbb T^3} |\vv_q|^2 - |\vell|^2 \dd x=
    \\&\int_{\mathbb T^3} \glue_i (|\vex_i|^2-|\vell|^2)+(1-\glue_i) (|\vex_{i+1}|^2-|\vell|^2) -2\glue_i(1-\glue_i)|\vex_i-\vex_{i+1}|^2 \dd x.
\end{align*}
The last term is $\lesssim \delta_{q+1}\ell_q^{2\alpha} $ by \eqref{e:vv_q-stability-L2}. As  $|\glue_i|\le1 $ and $ |1-\glue_i|\le 1 $, to prove \eqref{e:energy-of-vv_q}, we just need to show that for all $i=0,1,2,\dots,\imax $ and all $t\in I_{i-1}\cup J_i\cup I_i\subset[t_i-\tau_q,t_i+\tau_q]$,
\[ \left|\int_{\mathbb T^3} |\vex_i(x,t)|^2 - |\vell(x,t)|^2 \dd x\right| \lesssim \delta_{q+1}\ell_q^{2\alpha}. \]
For this, observe that $ \frac{\dd}{\dd t}\int_{\mathbb T^3}|\vell|^2 \dd x =  \frac{\dd}{\dd t}\int_{\mathbb T^3}|\vell|^2 - |\vex_i|^2 \dd x$, and\footnote{$A:B\coloneq A_{ij}B_{ij} $ is the double-contraction.}
\[ \frac{\dd}{\dd t}\int_{\mathbb T^3}|\vell|^2 \dd x 
=-2\int_{\mathbb T^3}\nabla  \vell : \Rell\dd x 
\lesssim  \|\vell\|_{H^1}  \|\Rell\|_\sobolev{3/2}{1} \ll 
  \tau_q^{-1} \delta_{q+1}\ell_q^{2\alpha},\] by the estimates in Proposition \ref{p:estimates-for-mollified}. Integrating from $t_i$ to $t$ gives
\[
\left|\int_{\mathbb T^3}|\vell(x,t)|^2 -|\vex_i(x,t)|^2 \dd x\right| \le  \tau_q\left| \frac{\dd}{\dd t}\int_{\mathbb T^3}|\vell(x,t)|^2 \dd x \right| \ll \delta_{q+1}\ell_q^{2\alpha},
\]
as claimed. Hence, the estimate \eqref{e:energy-of-vv_q} follows.
\end{proof}

\section{Definition of velocity increment $w_{q+1}$}
\label{s:wq+1-defn}

\subsection{Mikado flows}
\label{ss:mikado}
\newcommand{\lambA}{\lambda_+}
\newcommand{\lambB}{\lambda_- }
We first recall the  `geometric lemma', which was essentially first proven in \cite{zbMATH06210493}. For a quick proof using the inverse function theorem, see \cite{https://doi.org/10.48550/arxiv.1907.10436}.
\begin{lem}
There exist two disjoint sets of 6 directions $\Lambda_0,\Lambda_1\subset \mathbb Q^3\cap \sqrt 2 \mathbb S^2$ and a constant $\epsilon_\gamma\in(0,1)$ such that for each $k\in \Lambda_i$, $i=0,1$, there exist smooth  functions $\gamma_k : \mathbb B_{\epsilon_\gamma}(\Id) \to \mathbb R$ such that the following decomposition holds for all $R\in \mathbb B_{\epsilon_\gamma}(\Id)$:
\begin{align}
     R = \sum_{k\in\Lambda_i} \gamma_k^2(R)k\otimes k. \label{e:R-decomp}
\end{align}
In addition, we may extend each $k\in\Lambda_0\cup\Lambda_1$ to an orthogonal basis $(k,\kk,\kkk)$ of vectors in $\mathbb Q^3\cap \sqrt 2 \mathbb S^2 $so that no two of the 36 vectors are parallel.
\end{lem}

%
Next, we recall the Mikado flow construction as in e.g. \cite{zbMATH07301372}.
There exist points $x_k\in\mathbb T^3$ $(k\in \Lambda_1\cup\Lambda_2)$ and a constant $\epsilon_{\Lambda}<1/2$ such that the $\epsilon_{\Lambda}$-neighbourhoods of the periodic lines $L_k= \{x_k + tk\}_{t\in\mathbb R}\subset \mathbb T^3$ are pairwise disjoint.
Let $\Phi:\mathbb R^2\to \mathbb R$ be a smooth function compactly supported in $\mathbb B^{\mathbb R^2}_{\epsilon_\Lambda}(0),$  and such that $\phi\coloneq -\Delta \Phi$ solves
\[ \int_{\mathbb R^2}\phi(x_1,x_2)^2 \dd x_1 \dd x_2 = 1. \]
Also define  $\phi_{\lambB}: \mathbb R^2 \to \mathbb R$  using the $L^2(\mathbb R^2)$-invariant scaling, so that $\Delta \Phi_{\lambB} = \lambB^2  \phi_{\lambB}$:
\[ \phi_{\lambB}(x) \coloneq \lambB \phi(\lambB x), \quad \Phi_{\lambB}(x) \coloneq \lambB \Phi(\lambB x).\]
We identify $\phi_{\lambB}$ and $\Phi_{\lambB}$ with their $\mathbb T^2$-periodic extensions ($\mathbb T=\mathbb R/\mathbb Z$), supported on $\mathbb B^{\mathbb T^2}_{\epsilon_\Lambda/r}(0)$. We choose a sufficiently large integer $N_\Lambda$ so that  $N_{\Lambda}k, N_{\Lambda}\kk , N_{\Lambda}\kkk \in\mathbb Z^3$, and define $\phi_\bk,\Phi_\bk:\mathbb R^3\to \mathbb R$ by
\begin{align*}
    \phi_\bk(x) \coloneq \phi_{k,\lambB,\lambA }(x)\coloneq \phi_{\lambB}(\lambA N_{\Lambda} (x-x_k)\cdot \kk , \lambA N_{\Lambda}  (x-x_k)\cdot \kkk ),
    \\
    \Phi_\bk(x) \coloneq \Phi_{k,\lambB,\lambA }(x)\coloneq \Phi_{\lambB}(\lambA N_{\Lambda} (x-x_k)\cdot \kk , \lambA N_{\Lambda}  (x-x_k)\cdot \kkk ),
\end{align*}
which are supported on the $(\epsilon_{\Lambda}/r)$-neighbourhood of $L_k$, and $\frac{-1}{N_\Lambda^2 \lambA^2\lambB^2} \Delta \Phi_\bk = \phi_\bk$. Also, $\phi_\bk$ is $(\mathbb T/\lambA )^2$-periodic. Then,
\[ W_\bk(x)\coloneq W_{k,r,\lambA,\lambB}(x) \coloneq \phi_\bk (x) k \]
defines a family of $(\mathbb T/\lambA )^2$-periodic steady presureless Euler flows,  supported on the disjoint tubes $\{ d(x,L_k)<\epsilon_{\Lambda}/r\}$:
\begin{align}
    \begin{gathered}
    \div W_\bk = 0, \quad 
    \div (W_\bk \otimes W_\bk )  = 0, \quad 
    \int_{\mathbb R^3} W_\bk \dd x  = 0, \\
    \int_{\mathbb R^3} W_\bk \otimes W_\bk \dd x  = k\otimes k.
\end{gathered}\label{e:wbk-properties}
\end{align}
 hence, their linear combinations are also exact Euler flows. 
 
 
Also,\footnote{since $\curl (\nabla \Phi_\bk\times k)=\curl\curl( \Phi_\bk k)=\Delta  \Phi_\bk k- \nabla(k\cdot\nabla  \Phi_\bk)$ and $k\cdot\nabla  \Phi_\bk=0$.} $W_\bk = \frac1{\lambA\lambB} \curl V_\bk$ for the  vector field  $V_\bk$ supported on the tube $\{x: d(x,L_k)<\epsilon_\Lambda/r \}$,
\begin{align}
     V_\bk = \frac1{N_{\Lambda}^2\lambA\lambB}\nabla \Phi_\bk \times k.\notag 
\end{align}

We choose for the remainder of the paper
\begin{align}
    \lambA = \lambda_{q+1}^{5/6} , \quad \lambB = \lambda_{q+1}^{1/6},\quad   b\in 6\mathbb N. \label{e:constraint-intermittency} 
\end{align}
\begin{lem}[Integrability]\label{p:mikado-integrability}
For all $p\in[1,\infty]$, $N\ge 0$,
    \begin{align}
        \|\nabla^N W_\bk \|_{L^p(\mathbb T^3)} +   \|\nabla^N V_\bk \|_{L^p(\mathbb T^3)} \lesssim  \lambda  _{q+1}^{\frac16(1-\frac2 p)+N},    \end{align}
    where the implicit constant depends only on $ N_\Lambda, \epsilon_\Lambda,N,p,$ and $\phi$.
\end{lem}

\subsection{Space-time cutoffs}

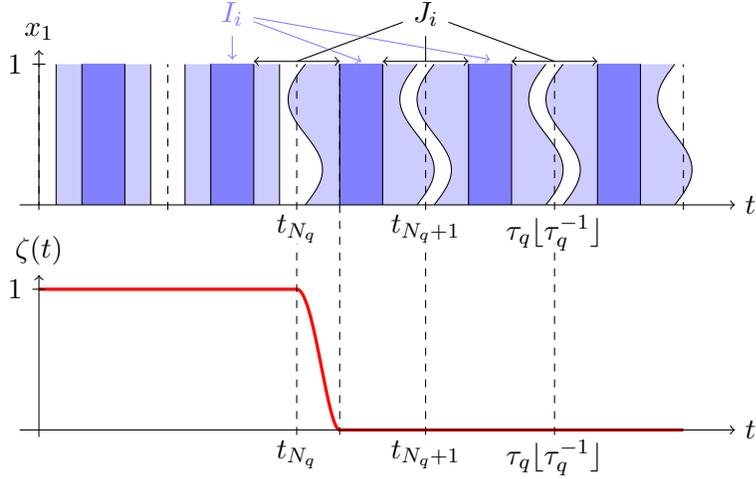
\begin{figure}
    \begin{tikzpicture}

\begin{axis}[clip=false,width=12cm,height=7.5cm,
    axis line style={draw=none}, 
    disabledatascaling,
    yticklabels={,,},     
    xticklabels={,,},
    tick style={draw=none},
    ]
 \addplot [line width = 0.4cm,domain=0:2*pi,samples=100](
	{\crvi{x}},
	{\crvii{x}}
	); 
 \addplot [white, line width = 0.35cm,domain=0:2*pi,samples=100](
	{\crvi{x}},
	{\crvii{x}}
	); 
 
\foreach \k [count=\K] in {0,...,\varN}
{
\ifnum \k<\varNi 
\addplot [domain=0:1,samples=100, name path global =A](%
    {\k*\varTAU+\crvi{x}},%
	{x}%
	)%
	node[%
        pos=0.5%
        ]%
        { }%
        ;
\addplot [domain=0:1,samples=100, name path global =B](
    {\k*\varTAU+\crvii{x}},
	{x}
	)
	node[%
        pos=0.5%
        ]%
        {}%
        ; 
\addplot [\colorEta] fill between [of=A and B];
\else 
\addplot [domain=0:1,samples=100, name path global = A](
    {\k*\varTAU+\crviii{x}},
	{x}
	)
	node[%
        pos=0.5%
        ]%
        { }%
        ;
        
\addplot [domain=0:1,samples=100, name path global = B](
    {\k*\varTAU+\crviv{x}},
	{x}
	)
	node[%
        pos=0.5%
        ]%
        {}%
        ; 
\addplot [\colorEta] fill between [of=A and B];
\fi 
\addplot [domain=-0.035:1,samples=2, dashed]( 
    {\k*\varTAU},
	{x}
	)   ;
\addplot [domain=0:1,samples=2, name path global = A](
    {\k*\varTAU+\varTAU/3},
	{x}
	)
        ;
\addplot [domain=0:1,samples=2, name path global = B](
    {\k*\varTAU+2*\varTAU/3},
	{x}
	)
	node[%
        pos=0.5%
        ]%
        { }%
        ;
\addplot [\colorR] fill between [of=A and B];
}
\addplot [domain=-0.035:1,samples=2, dashed]( 
    {\varN*\varTAU+\varTAU},
	{x}
	)
	node[%
        pos=0.5%
        ]%
        { }%
        ;

\draw[->] (-0.15*\varTAU,0) -- (\varN*\varTAU+1.4*\varTAU,0) node[right] {$t$};
\node at (\varNi*\varTAU,0) [below]{$t_{N_q}$};
\node at (\varNi*\varTAU+ \varTAU,0) [below]{$t_{N_q+1}$};
\node at (\varNi*\varTAU+2*\varTAU,0) [below]{$\tau_q\lfloor\tau_q^{-1}\rfloor$};

\draw [->] (0,0) -- (0,1+\varTAU/2) node [above] {$x_1$} ;
\draw (\varTAU/20,1) -- (-\varTAU/20,1) node [left] {$1$} ;

\node at (\varNi*\varTAU+\varTAU/3+\varTAU/6-\varTAU,1+\varTAU) (a) [above, \colorR]{$I_i$};
\draw [\colorR,->] (a) -- (\varTAU*2.5-\varTAU,1+\varTAU/6);
\draw [\colorR,->] (a) -- (\varTAU*3.5-\varTAU,1+\varTAU/10);
\draw [\colorR,->] (a) -- (\varTAU*4.5-\varTAU,1+\varTAU/6); 

\node at (\varNi*\varTAU+\varTAU*2-\varTAU,1+\varTAU ) (b) [above]{$J_i$};
\draw [thick, white] (b) -- (\varNi*\varTAU,1+\varTAU/10) node (c) {} ;
\draw [] (b) -- (\varNi*\varTAU,1+\varTAU/10) node (c) {} ;
\draw [<->] ($(c)-(\varTAU/3,0)$) -- ($(c)+(\varTAU/3,0)$);
\draw [thick, white] (b) -- (\varNi*\varTAU+2*\varTAU-\varTAU,1+\varTAU/10) node (c) {} ;
\draw [] (b) -- (\varNi*\varTAU+2*\varTAU-\varTAU,1+\varTAU/10) node (c) {} ;
\draw [<->] ($(c)-(\varTAU/3,0)$) -- ($(c)+(\varTAU/3,0)$);
\draw [thick, white] (b) -- (\varNi*\varTAU+3*\varTAU-\varTAU,1+\varTAU/10) node (c) {} ;
\draw [] (b) -- (\varNi*\varTAU+3*\varTAU-\varTAU,1+\varTAU/10) node (c) {} ;
\draw [<->] ($(c)-(\varTAU/3,0)$) -- ($(c)+(\varTAU/3,0)$);


\draw [\colorEtaMinusOne, very thick] (0,\varYOFF+1) -- (\varETASTART +0.0005,\varYOFF+1);
\addplot [\colorEtaMinusOne, very thick, domain=\varETASTART :(\varETASTART +\varTAU/3),samples=100] ( 
    {x},
	{\etaminusone{(x-(\varETASTART) )/(\varTAU/3) } +\varYOFF}
	)
	node[%
        pos=0.5%
        ]%
        { }%
        ;
\draw [\colorEtaMinusOne, very thick] (\varETASTART +\varTAU/3-0.0005,\varYOFF) -- (\varN*\varTAU+\varTAU,\varYOFF)  ;
\draw[->] (-0.15*\varTAU,\varYOFF ) -- (\varN*\varTAU+1.4*\varTAU,\varYOFF ) node[right] {$t$};
\node at (\varNi*\varTAU, \varYOFF ) [below]{$t_{N_q}$};
\node at (\varNi*\varTAU+\varTAU,\varYOFF ) [below]{$t_{N_q+1}$};
\node at (\varNi*\varTAU+2*\varTAU,\varYOFF ) [below]{$\tau_q\lfloor\tau_q^{-1}\rfloor$};

\draw [->] (0,\varYOFF -0.05) -- (0,\varYOFF + 1+\varTAU/2 ) node [above] {$\zeta(t)$} ;

\addplot [domain=-0.035:1-0.2*\varYOFF,samples=2, dashed]( 
    {\varNi*\varTAU},
	{x+\varYOFF}
	)
	node[%
        pos=0.5%
        ]%
        { }%
        ;
\addplot [domain=-0.035:1-1*\varYOFF,samples=2, dashed]( 
    {\varNi*\varTAU+\varTAU/3},
	{x+\varYOFF}
	)
	node[%
        pos=0.5%
        ]%
        { }%
        ;

\addplot [domain=-0.035:1-0.2*\varYOFF ,samples=2, dashed]( 
    {\varNi*\varTAU+\varTAU},
	{x+\varYOFF}
	)
	node[%
        pos=0.5%
        ]%
        { }%
        ; 
\addplot [domain=-0.035:1-0.2*\varYOFF ,samples=2, dashed]( 
    {\varNi*\varTAU+2*\varTAU},
	{x+\varYOFF}
	)
	node[%
        pos=0.5%
        ]%
        { }%
        ;
\draw (\varTAU/20,1+\varYOFF) -- (-\varTAU/20,1+\varYOFF) node [left] {$1$} ;
\end{axis} 
\end{tikzpicture}
    \caption{The support of the space-time cutoffs.  $\zeta$ is used in \eqref{e:rho-q-i} to transition between the two types of cutoffs. This figure is taken from \cite{https://doi.org/10.48550/arxiv.2204.03344}.}
        \label{f:space-time}
\end{figure}

The point of this subsection is to formalise Figure \ref{f:space-time}. These transition from Isett's cutoffs \cite{zbMATH06976275} to the $S$-shaped cutoffs of Buckmaster--Vicol \cite{zbMATH07003146}, and was already used in \cite{https://doi.org/10.48550/arxiv.2204.03344}.

\label{ss:spacetime-cutoffs}
 Define the index
$ i_{q} \coloneq  \big\lfloor \frac1{\tau_q} \big\rfloor - 2\in\mathbb Z_{\ge 0}.  $
Let $\bar\eta_0\in C_c^\infty(J_0\cup I_0 \cup J_1;[0,1])$  satisfy
$  \supp \bar\eta_0 = I_0 + \big[  {-\frac{\tau_q}6},\ \frac{\tau_q}6\big] $
be identically 1 on $I_0$, and satisfy the derivative estimates for $N\ge 0$:
\[ \|\del_t ^N \bar\eta_0\|_{C^0_t} \lesssim \tau_q^{-N}. \]
Then, we set $\bar\eta_i(t)\coloneq\bar\eta_0(t-t_i)$ for $0\le i\le i_{q}$. 
For $i_q\le i\le \imax$,
let $\epsilon \in (0,1/3)$,  $\epsilon_0\ll 1$, fixed throughout the iteration scheme, and define  
\begin{align*}
     I_{i}'&\coloneq  \Big[t_i+\frac{\epsilon\tau_q}3, t_i+ \frac{(3-\epsilon)\tau_q}3\Big]  \subset\mathbb R ,\\
         S_i &\coloneq \Big\{\big(x,t+\frac{2\epsilon\tau_q}3\sin (2\pi x_1) \big) : x\in\mathbb T^3, \ t \in I_i'\Big\}\subset  \mathbb T^3 \times \mathbb R,   \\
         \tilde \eta_i (x,t)&\coloneq \frac1{\epsilon_0^3}\frac1{\epsilon_0 \tau_q} \iint_{S_i}\phi\Big(\frac{|x-y|}{\epsilon_0}\Big)\phi\Big(\frac{|t-s|}{\epsilon_0\tau_q}\Big)\dd y \dd s.
\end{align*}
We define the space-time cutoffs $\eta_i$, $i\ge 0$, by
\begin{align}
    \eta_i (x,t) \coloneq \begin{cases}
 \bar\eta_i(t), & 0\le i< i_{q}, \\
 \tilde\eta_i(x,t), & \phantom{0\le{}}i\ge i_{q}.
 \end{cases} \label{e:eta-defn}
\end{align}

\begin{lem}\label{p:eta}
 The  functions $\eta_i$ satisfies  for all $i\le \imax$ (see Figure \ref{f:space-time}):
    \begin{enumerate}
    \item $\eta_i \in C^\infty_c(\mathbb T^3 \times (J_i\cup I_i\cup J_{i+1}) ; [0,1])$, and for $n,m\ge 0$:\begin{align}
         \|\del^n_t  \eta_i\|_{C^m   }  \lesssim_{n,m} \tau_q^{-n}, \label{e:eta-estimates}
    \end{align}
    \item $\eta_i(\cdot,t) \equiv 1$ for $t\in I_i$,
    \item $\supp \eta_i$ are pairwise disjoint,
    \item  for all $t\in[t_{i_{q}},T]$, $c_\eta:=\frac14  \le \sum_{i=0}^{\imax}\int_{\mathbb T^3}  \eta_i^2(x,t) \dd x \le 1$,
    \item for all $0\le     i<i_{q}$, $\eta_i $ does not depend on $x$, and  is supported in time on the $\tfrac{\tau_q}6$-neighbourhood of $I_i$.
\end{enumerate}
\end{lem}

\subsection{Stress cutoffs}
%
Here, we recall the stress cutoffs from \cite{zbMATH07003146} to obtain $L^\infty$ control on $\RRR_q$. Let $0\le \tilde \chi_0,\tilde \chi \le 1$ with $\tilde\chi_0$ supported on $[0,4]$, $\chi_0\equiv1$ for $x\le 1/4$, $\tilde\chi^2 = \tilde \chi_0(x)^2 - \tilde \chi_0(4x)^2$.  Then their rescalings give a partition of unity:
\begin{align}
 \widetilde{\chi}_{0}^{2}(y)+\sum_{j \geq 1} \tilde{\chi}_{j}^{2}(y) \equiv 1, \quad \text { where }  \quad\widetilde{\chi}_{j}(y)\coloneq \widetilde{\chi}\left(4^{-j} y\right), \quad  j\ge 1, \, y\ge0.\notag 
 \end{align}

Using the `Japanese bracket' $\left\langle v\right\rangle \coloneq   \sqrt{1+|v|^2}$), we define 
\begin{align}
\chi_\bj(x, t)\coloneq \chi_{j, q+1}(x, t)\coloneq  \widetilde{\chi}_{j}\Bigg(\Bigg\langle\frac{\RRR_q (x, t)}{\yibai\lambda_q^{-\alpha}\delta_{q+1}}\Bigg\rangle\Bigg),\notag 
\end{align}
which is also a partition of unity, i.e. for all $x\in\mathbb T^3,t\ge0$,
\[ \sum_{j\ge 0} \chi^2_\bj (x,t)\equiv 1.\]

\begin{prop}On the support of $\chi_\bj$, $4^{j-1}-1  \le  \frac{|\RRR_q|}{\yibai \lambda _q^{-\alpha}\delta_{q+1}} \le 4^{j+1}$. Hence, there is a final $\jmax=\jmax(q) $ such that $\chi_\bj\equiv 0$ for all $j\ge \jmax$. In addition, for all $p\in[1,\infty]$,\label{p:chi}
\begin{gather}
4^{\jmax} \lesssim  \ell_q^{-3}, \notag 
\\
    \|\chi_\bj \|_{L^1} \ll 4^{-j}, \quad     \|\chi_\bj \|_{L^2} \ll 2^{-j}, 
    \qquad \| \chi_\bj \|_{L^p} \lesssim 4^{-j/p }, \label{e:chi-lp}\\ 
    \phantom{\del_t}\|\nabla^n\chi_\bj \|_{L^p} \lesssim  4^{-j(N+1/p)}\ell_q^{-4N},
    \phantom{\tau_q^{-1}\ell_q^{-3}}
     \qquad  (n\ge 1)  , \label{e:chi-Dn-Lp} 
     \\
    \|\nabla^n \partial_t\chi_\bj \|_{L^p} \lesssim   4^{-j(N+1/p)}\ell_q^{-4N}\tau_q^{-1} \ell_q^{-3},\qquad (n\ge 1). \label{e:chi-DnDt-Lp}
\end{gather}

\end{prop}
\begin{proof} The embedding $\sobolev{3+\alpha/10 }{1}\subset L^\infty $ and \eqref{e:RRR_q-WN1} give the upper bound  \[ \|\RRR_q \|_{L^{\infty}} \le \|\RRR_q\|_\sobolev {3+\alpha/10}1\le \delta_{q+1}\ell_q^{-3-\alpha/10}\lambda_q^{-2\alpha} \ll \delta_{q+1}\ell_q^{-3}\lambda_q^{-\alpha}.\] 
 From the  lower bound of the support, it suffices to take as $\jmax$ the least $j$ such that $4^{j-1}\ge  \ell_q^{-3} $.
 It follows that $4^{\jmax}\le 2  \ell_q^{-3}$.
  
The $L^p$ estimates \eqref{e:chi-lp} follow by interpolating between $\|\chi\|_{L^\infty}\le1$ and the  bound
\[ |\supp \chi_\bj| \le \Bigg|\Bigg\{ x:\frac{|\RRR_q|}{\yibai\lambda_q^{-\alpha} \delta_{q+1}}> 4^j  \Bigg\}\Bigg| \le\frac{ \delta_{q+1} \lambda_q^{-2\alpha} }{ \yibai 4^j \lambda_q^{-\alpha}\delta_{q+1} } \ll 4^{-j} .\]
\cite[Prop C.1]{zbMATH06456007} and the embedding $W^{3,1+\alpha}\subset L^\infty$ give
\begin{multline*}
     \|\nabla^N \chi_\bj \|_{L^\infty} \lesssim 4^{-j} \|\nabla^N\RRR_q\|_{L^\infty} \lambda_q^{\alpha}\delta_{q+1}^{-1} + 4^{-Nj}\|\nabla \RRR_q\|^N_{C^0} \lambda_q^{N \alpha}\delta_{q+1}^{-N}.
     \\
     \lesssim 4^{-Nj}(\delta_{q+1}\ell_q^{-4+\alpha} )^N\lambda_q^{N\alpha} \delta_{q+1}^{-N} \lesssim 4^{-Nj}\ell_q^{-4 N}(\ell_q\lambda_q)^{N\alpha} \ll 4^{-Nj}\ell_q^{-4 N}.
\end{multline*}
Interpolating again with the support bound gives an extra power of $4^{-j}$. A similar calculation shows the inequality for $\del_t\chi_\bj$:
\begin{align*}
\MoveEqLeft    \|\nabla^N \del_t \chi_\bj\|_{L^\infty} \lesssim  \lambda_q^{\alpha}\delta_{q+1}^{-1}\left\| \nabla^N\bigg(\tilde\chi_j(\langle\bullet\rangle)'\circ \bigg(\frac{\RRR_q}{\yibai\lambda_q^{-\alpha}\delta_{q+1}}\bigg) \del_t\RRR_q\bigg)\right\|_{L^\infty}
    \\
    &\lesssim \lambda_q^{\alpha}\delta_{q+1}^{-1}\Big( \|\del_t\RRR_q\|_{L^\infty}(4^{-2j}\|\nabla^N \RRR_q\|_{L^\infty} +4^{-(N+1)j}\|\nabla\RRR_q\|_{C^0}^N\lambda_q^{N\alpha}\delta_{q+1}^{-N})\\
    &\quad  +4^{-j} \|\del_t\nabla^N\RRR_q\|_{L^\infty}\Big)
    \\
    &\lesssim \lambda_q^{\alpha}\delta_{q+1}^{-1}\Big( \tau_q^{-1}\delta_{q+1}\ell_q^{-3+\alpha}\Big(4^{-2j}\delta_{q+1}\ell_q^{-3-N+\alpha} 
    \\
    &\quad +4^{-(N+1)j}(\delta_{q+1}\ell_q^{-4+\alpha})^N\lambda_q^{N\alpha}\delta_{q+1}^{-N}\Big)  + 4^{-j}\delta_{q+1} \tau_q^{-1}\ell_q^{-N+\alpha}\Big)
    \\
    &\lesssim 4^{-(N+1)j} \tau_q^{-1}\ell_q^{-3-4N} (\ell_q\lambda_q)^{2\alpha}
    \\
    & \ll  (4^{-(N+1)j} \tau_q^{-1}\ell_q^{-3})\ell_q^{-4N}.
\end{align*}
We used that the term with $4^{-(N+1)j}$ is largest because $\ell_q^{-3}4^{-j}>1$.
%
\end{proof}
\begin{lem} \label{l:chi-0-integral}$\int_{\mathbb T^3} \chi_{(0)}^2 \dd x = \int_{\mathbb T^3} \widetilde{\chi}_{j}\Big(\Big\langle\frac{\RRR_q (x, t)}{\yibai\lambda_q^{-\alpha}\delta_{q+1}}\Big\rangle\Big) \dd x \ge 1/2$ for $a\gg 1$.
\end{lem}
\begin{proof} Write $\mathsf R=\frac{|\RRR_q|}{\yibai\lambda_q^{-\alpha} \delta_{q+1}}$ for brevity. By Markov's inequality, for any $c>0$,
    \begin{align*}
        \int_{\mathbb T^3} \chi_{(0)}^2 \dd x 
        &\ge \int_{ |\mathsf R| < c } \tilde\chi_0(\langle \mathsf R\rangle ) ^2\dd x 
        \\
        &\ge \inf_{s\in [0,\langle c \rangle )} \tilde\chi_0^2(s) \lvert\{ |\mathsf R | < c \}| 
        \\
        &= \inf_{s\in [0,\langle c \rangle )} \tilde\chi_0^2(s)  \Big(1 - |\{ |\mathsf R | \ge c \}|\Big) 
        \ge  \inf_{s\in [0,\langle c \rangle )} \tilde\chi_0^2(s)\bigg(1- \frac{\|\mathsf R\|_{L^1}}{c}\bigg).
    \end{align*}
    Since $\|\mathsf R\|_{L^1} \lesssim  \lambda_q^{-\alpha} \to 0 $ as $a\to\infty$, choosing first $c\ll 1$ so that $\inf_{[0,\langle c \rangle )}\tilde\chi_0^2\ge 2/3$, we can take $a$ large to obtain $1/2$ as a lower bound.
\end{proof}

\subsection{Energy gap }
\label{ss:energy-gap}
Define for $j\ge1$,\begin{align*}
\rho_{q,j}\coloneq 4^{j+j_\gamma}\delta_{q+1}\lambda_q^{-\alpha},
\end{align*}
where $j_\gamma$ is a fixed geometrical constant chosen so that 
\begin{align}
    \bigg \| \frac{\RRR_q}{\rho_{q,j}}\bigg \|_{L^\infty(\supp \chi_\bj)} \le \epsilon_\gamma . \label{e:TR-j-welldef-jge1}
\end{align}
\eqref{e:TR-j-welldef-jge1} ensures that $\gamma_k(\TR{j})$ is well-defined in Subsection \ref{sss:principal-and-cor} below. From Proposition \ref{p:chi}, 
$   \| \RRR_q \|_{L^\infty(\supp \chi_\bj)} \le 4^j \lambda_{q}^{-\alpha}  \delta_{q+1} = \rho_{q,j} 4^{-j_\gamma}. 
$ Hence, one can take $j_\gamma = \log_4\frac1{\epsilon_\gamma}$  which is independent of all parameters $b,\beta,\alpha,a$.
 
   For $j=0$, we set
\begin{align}
\tilde\rho_{q,0}(t) &\coloneq   \frac{1}3\Big(e(t) - \int_{\mathbb T^3} |\vv_q|^2 \dd x-3\sum_{j\geq 1,i} \rho_{q,j} \int_{\mathbb T^3}\eta^2_{ i}\chi_\bj^2 \dd x - \frac{\delta_{q+2}}2 \Big) .\label{e:rho-tilde-q0-defn}
\end{align}
As  $\eta_{i }|_{\supp \RRR_q}\equiv  1$, $4^j\lesssim \langle \frac{\RRR_q}{\yibai\lambda_q^{-\alpha}\delta_{q+1}}\rangle$ on the support of $\chi_\bj$, and $\sum_j \chi_\bj^2 \equiv 1$,
\[ \sum_{j=1}^{\jmax }  4^j \int_{\mathbb T^3} \chi^2_\bj \dd x \le \int 1+ \frac{|\RRR_q|}{\yibai\lambda_q^{-\alpha}\delta_{q+1}} \dd x \lesssim  1.  \]
It follows that $\sum_{j\geq 1,i} \rho_{q,j} \int_{\mathbb T^3}\eta^2_{ i}\chi_\bj^2 \dd x \lesssim \delta_{q+1}\lambda_{q}^{-\alpha}$ uniformly in $t$, so that with the inductive estimate \eqref{e:energy-vell} and Proposition \ref{p:energy-of-vv_q}, we obtain (for $a\gg 1$)
\begin{align}
\frac{1}{2}\delta_{q+1}\le  \tilde\rho_{q,0}(t) \le  2\delta_{q+1},\quad\forall t\in [\tau_{q-1},T]. \label{e:rho-tilde-q-0-bd}
\end{align}
We further define a function 
$  \zeta\in C_c^\infty([0,t_{i_{q}+1});[   0,1])$ as in Figure \ref{f:space-time}
 such that 
$
 \zeta \equiv 1    
$ 
 for $ t\le t_{i_{q}}$, 
 and it satisfies the estimates for all $N\ge 0$,
\begin{align*}
   \sup_t |\del_t^N \zeta(t)| \lesssim \tau_q^{-N}. 
\end{align*}
Note that $t_{i_{q}+1}=\tau_q (\lfloor \tau_q^{-1}\rfloor -1) \le 1-\tau_q$, so that $t\ge 1-\tau_q$ implies $\zeta(t)=0$.
We decompose $\tilde\rho_{q,0}$ by setting 
\begin{align}
 \rho_{q,0} (t) &\coloneq \delta_{q+1}\zeta(t)+\frac{ \tilde\rho_{q,0}(t)(1-\zeta(t))}{ \zeta(t) + \sum_{i= i_{q}}^{\imax} \int_{\mathbb T^3}\eta_i^2\chi_{(0)}^2 (\tilde x,t)\dd \tilde x },  \label{e:rho-q-i}
\end{align}
Observe that $\rho_{q,0}$  is well-defined: indeed,
 when
$t\in I_i$ with $i\geq i_q$, we have the lower bound $\int_{\mathbb T^3}\eta_i^2\chi_{(0)}^2 (\tilde x,t)\dd \tilde x=\int_{\mathbb T^3}\chi_{(0)}^2 (\tilde x,t)\dd \tilde x\ge \frac12$, by \eqref{l:chi-0-integral}.
 If instead $t\in J_i$ with $i\geq i_q$, then $\RRR_q=0$ and $\chi_{(0)}=1$, so $\int_{\mathbb T^3}\eta_i^2\chi_{(0)}^2 (\tilde x,t)\dd \tilde x=\int_{\mathbb T^3}\eta_i^2 (\tilde x,t)\dd \tilde x\ge c_{\eta}= \frac{1}{4} $. Therefore,
\begin{align}
    \zeta(t) + \sum_{i= 0}^{\imax} \int_{\mathbb T^3}\eta_i^2\chi_{(0)}^2 (\tilde x,t)\dd \tilde x  \geq \begin{cases}
 1 ,& 0\le t\le  t_{i_{q}}, \\
1/4 ,& t_{i_{q}}\le  t\le  T,
 \end{cases} \label{e:rho-q0-denominator-lowerbd}
\end{align}
and the fraction in \eqref{e:rho-q-i} makes sense. 

\begin{prop}\label{p:rho-estimates}
 For all  $t \in[0,T]$, and each $p\in[1,\infty)$,
    \begin{gather} \frac{\delta_{q+1}}{4} \le \rho_{q,0}(t)  \le  8\delta_{q+1} ,\label{e:rho-q0-upper-lowerbd} \\
    \|\partial_{t} \rho_{q,0}\|_{L^\infty_t}  \lesssim \tau^{-1}_q\delta_{q+1}, \ \text{and}\label{e:rho-q0-Dt}\\
    \sum_{j=0}^{\jmax} \|\rho_{q,j}\|_{L^\infty_t}^{1/p} 4^{-j/p  } \le 9^{1/p} \delta_{q+1}^{1/p}.\label{e:rho-sum-bd}
    \end{gather}
In particular,  \eqref{e:rho-q0-upper-lowerbd} implies that for $a\gg 1$,  $\Big \| \frac{\RRR_q}{\rho_{q,0}}\Big\|_{L^\infty(\supp \chi_{(0)})} \le \epsilon_\gamma $, so that $\gamma_k(\TR 0)$ is well defined in Subsection \ref{sss:principal-and-cor} below.  
\end{prop}  
\begin{proof}
Since $\zeta(t) + \sum_{i= 0}^{\imax} \int_{\mathbb T^3}\eta_i^2\chi_{(0)}^2 (\tilde x,t)\dd \tilde x  \le  2$, we have by \eqref{e:rho-tilde-q-0-bd} that
\begin{align}
    \rho_{q,0}(t)  \geq \begin{cases}
 \delta_{q+1}, & 0\le t\le  1-\tau_{q}, \\
 {\delta_{q+1}}/{4}, &1- \tau_{q}\le  t\le  T.
 \end{cases} \label{e:rho-q0-lowerbd}
\end{align}
Also, $\rho_{q,0}(t)\le  8\delta_{q+1}$ by \eqref{e:rho-tilde-q-0-bd} and \eqref{e:rho-q0-denominator-lowerbd}, which proves \eqref{e:rho-q0-upper-lowerbd}. 
 \eqref{e:rho-q0-Dt} follows by a direct calculation with Proposition \ref{p:chi}. For \eqref{e:rho-sum-bd}, with $p\in[1,\infty)$ fixed,
\begin{align}
    \sum_{j=0}^{\jmax} \|\rho_{q,0}\|_{L^\infty_t}^{1/p} 4^{-j/p}
    &= \|\rho_{q,j}\|_{L^\infty_t}^{1/p} + \sum_{i=1}^{\jmax} \rho_{q,j}^{1/p} 4^{-j/p}\notag\\
    &\le 8^{1/p} \delta_{q+1}^{1/p}+\jmax 4^{j_\gamma/p} \delta_{q+1}^{1/p}\lambda_{q}^{-\alpha/p}\notag\\
    &\le 8^{1/p} \delta_{q+1}^{1/p}+ 4^{j_\gamma/p} \delta_{q+1}^{1/p}\lambda_{q}^{-\alpha/p}\lvert\log_4 \ell_q^{3/2}| \le 9^{1/p} \delta_{q+1}^{1/p},\notag 
    \end{align}
    where the last inequality holds by choosing $a$ sufficiently large:
\[ \log_4\frac1{\ell_q} \le C \log_4\frac1{\lambda_q} \le C_{\alpha,p} \lambda^{\alpha/{2p}}_q,\] so we can use $C_{\alpha,p} \lambda_q^{\alpha/2p} \ll 1$ with $a\gg 1$.
\end{proof}

\subsection{Principal part and corrector terms}
\label{sss:principal-and-cor}
Define 
\begin{align}
          \TR j\coloneq \Id-\frac{\RRR_q(t,x)}{\rho_{q,j}(t)} ,\label{e:TRj-defn}\quad\text{and}\quad a_\bk \coloneq a_{j,i,k}(\TR j) \coloneq \eta_{ i}\chi_\bj\rho_{q,j}^{1/2}  \gamma_k(\TR j). 
 \end{align}
The \emph{principal part} of the perturbation is
\begin{align*}
    \wpq \coloneq \sum_{j=0}^{\jmax}\sum_{i=0}^{\imax}\sum_{k \in \Lambda_i} a_\bk  W_\bk.  
\end{align*}
As $a_\bk$ is not constant, $\wpq$ is not divergence-free (equivalently, not a curl). We fix this using  identity $\curl(a_\bk V_\bk )=\nabla a_\bk\times V_\bk +a_\bk \curl V_\bk $. That is, we define\begin{align}
    \wcq &\coloneq \frac1{\lambda_{q+1}} \sum_{j,i,k} \nabla a_\bk \times V_\bk . \notag
\end{align}
Then $\wcq+\wpq=\frac1{\lambda_{q+1}} \curl\left(\sum_{j,i,k} a_\bk V_\bk\right) $ is divergence free with zero mean.

Finally, we define the data corrector $\wdq $ by
\begin{align}
    \wdq  &\coloneq \vin * \psi_{\ell_q}-\vin * \psi_{\ell_{q-1}}* \psi_{\ell_q}. \notag
\end{align}
Note that $\wdq $ is divergence free and has zero mean. We define the perturbation $w_{q+1}$ and the new velocity $v_{q+1}$ by
\begin{align}
    w_{q+1} &\coloneq \wpq +\wcq+\wdq ,   \quad\text{and}\quad v_{q+1} \coloneq \vv_q + w_{q+1}. \notag
\end{align}

\begin{lem}[Estimates for the amplitude]\label{p:amplitude} Let $U_j=\supp \chi_\bj$.
For all $N\ge 1$ and all $p\in[1,\infty)$,
\begin{align}
\| \nabla^N \gamma_k(\TR  j)\|_{L^\infty (U_j)} &\lesssim \ell_q^{-4N}, \label{e:gamma-Dn-C0}
\\
\| \nabla^N\del_t \gamma_k(\TR  j)\|_{L^\infty (U_j)} &\lesssim \ell_q^{-4N}  \tau_q^{-1}\ell_q^{-3} \label{e:gamma-DnDt-C0},
\\
\|\nabla^N a_\bk\|_{L^p(U_j)} &\lesssim \ell_q^{-4N} \|\rho_{q,j}^{1/2}\|_{L^\infty_t} 4^{ -j/p} \lesssim \ell_q^{-4N}\delta_{q+1}^{1/2} 4^{ j(\frac12-\frac1p)},   \label{e:a-bk-Dn-Lp} 
\\
\|\nabla^N\del_t a_\bk \|_{L^p(U_j)} &\lesssim \ell_q^{-4N}\delta_{q+1}^{1/2} \tau_q^{-1}\ell_q^{-3} 4^{j(\frac12 -\frac1p)}. \notag
\end{align}

\end{lem}
\begin{proof}
The estimates for $ \gamma_k(\TR j)$ on the support of $\chi_\bj$ are proven similarly to the estimates \eqref{e:chi-Dn-Lp} and \eqref{e:chi-DnDt-Lp} of $\chi_\bj$.
Since
\begin{align*}
    \|\nabla^N( \eta_{ i}\chi_\bj\rho_{q,j}^{1/2}) \|_{L^p} 
    &\le 
    \|\rho_{q,j}^{1/2}\|_{L^\infty_t}  (\|\nabla^N\eta_i\|_{L^\infty} \|\chi_\bj\|_{L^p}+ \|\eta_i\|_{L^\infty} \|\nabla^N\chi_\bj\|_{L^p}),
\end{align*}
the estimate \eqref{e:a-bk-Dn-Lp} follows from \eqref{e:gamma-Dn-C0}, \eqref{e:eta-estimates}, and \eqref{e:chi-Dn-Lp}. 

Finally, the estimate for $\nabla^N\del_t a_\bk$ follows by similarly estimating each of the terms in
\begin{align*}
    \del_t a_\bk 
     &= (\del_t\eta_{ i})\chi_\bj\rho_{q,j}^{1/2}\gamma_k(\TR j)
    + \eta_{ i}(\del_t \chi_\bj)\rho_{q,j}^{1/2}\gamma_k(\TR j)
    \\
    &\quad + \eta_{ i}\chi_\bj(\del_t\rho_{q,j}^{1/2})\gamma_k(\TR j)
    + \eta_{ i}\chi_\bj\rho_{q,j}^{1/2}(\del_t\gamma_k(\TR j)),
\end{align*}
using \eqref{e:eta-estimates}, \eqref{e:chi-DnDt-Lp}, \eqref{e:rho-q0-Dt}, and \eqref{e:gamma-DnDt-C0}.
\end{proof}

\begin{cor}[Estimates for $\wpq$, $\wcq$ and $\wdq$] Let $M=361$.
\label{c:estimates-for-wpq-wcq}
\begin{align}
    \|\wpq\|_{L^2}+\frac1{\lambda_{q+1}} \|\wpq\|_{H^1} &\le \frac{M}{4}\delta_{q+1}^{1/2},\label{e:wpq-L2}\\
        \|\wcq\|_{L^2}+\frac1{\lambda_{q+1}} \|\wcq\|_{H^1} &\le \ell_q^{-4}\lambda_{q+1}^{-1} \delta_{q+1}^{1/2},\label{e:wcq-L2}\\
        \|\wdq \|_{L^2}+\frac1{\ell_q^{-1}} \|\wdq\|_{H^1} &\le 3\lambda_{q+1}^{-4\alpha} \delta_{q+2} \ll \delta_{q+1}^{1/2},\label{e:wdq-L2}\\
            \|w_{q+1}\|_{L^2}+\frac1{\lambda_{q+1}} \|w_{q+1} \|_{H^1} &\le \frac M2 \delta_{q+1}^{1/2}\label{e:wq+1-L2},\\
            \|\wpq \|_{L^4}+\|\wcq \|_{L^4}+\|\wdq \|_{L^4} &\lesssim  \delta_{q+1}^{1/2} \lambda_{q+1}^{\frac{3/4}6} + \frac{\lambda_{q+1}^{-1+\frac{3/4}6}}{\ell_q^4}  + \frac{\ell_{q-1}  }{\ell_q^2} \ll \lambda_{q+1}^{4}. \label{e:wq+1-L4-crude}
\end{align}
In particular, \eqref{e:v-q-H0} and \eqref{e:v-q-H1} hold with $q$ replaced with $q+1$, once $a\gg1$.
\end{cor}
\begin{proof}
\eqref{e:wpq-L2} and \eqref{e:wcq-L2} follow from the improved \Holder inequality (Lemma \ref{p:improved-holder-BV}), Lemma \ref{p:amplitude}, Lemma \ref{p:mikado-integrability}, \eqref{e:rho-sum-bd}, and enforcing $\ell_q^{-4} \ll \lambda_{q+1}$. Specifically for \eqref{e:wpq-L2}, as $|\Lambda_i|=6$,
\[ \|\wpq \|_{L^2}\le 5\sum_{ k\in \Lambda_i}\Bigg\lVert\sum_{\substack{j,i\ge0}}a_\bk \Bigg \|_{L^2} \|W_\bk\|_{L^2} \le 30 \sup_{i\ge 0}  \sum_j \|\rho_{q,j}\|_{L^\infty}^{1/2}2^{-j} \le 90\delta_{q+1}^{1/2}.  \]
The gradient has two terms, one about the same size (after dividing by $\lambda_{q+1}$) as $\|\wpq\|_{L^2}$ when the gradient falls on $W_\bk$, the other much smaller if it falls on $a_\bk$ (when $a\gg1$), so we can take $M=90\cdot4 +1=361$.  To ensure that Lemma \ref{p:improved-holder-BV} is applicable with parameters  $\mu = \ell_q^{-4}$, $\kappa = \lambda_{q+1}^{5/6}$, we need
$\ell_q^{-4} \le \frac1{\sqrt{27}} \lambda_{q+1}^{5/6}$. For this, it is enough to enforce
\begin{align}
    b>\frac{4-4\beta}{\frac56-4\beta}. \label{e:constraint-b}
\end{align}Then we can take $\alpha\ll1$,  $a\gg1$, and $K\gg1$. The parameter $K$ is from Lemma \ref{p:improved-holder-BV} and can depend on the parameters $a,b,q$.
For example, we can take 
$\beta<1/21$ and $b\ge 6$, which is a weaker constraint than \eqref{e:constraint-b-beta-strong} below.


For \eqref{e:wdq-L2}, since $\|\vin\|_{H^{\betain}}\le 1$, 
\begin{align*}
    \|\wdq \|_{L^2} &= \|\vin* \psi_{\ell_q}-(\vin* \psi_{\ell_{q}})* \psi_{\ell_{q-1}}\|_{L^2}
    \le   \ell^{\betain}_{q-1},
\end{align*}
and we need to show that this is controlled by $\lambda_{q+1}^{-4\alpha} \delta_{q+2}$.  If we just bound $\ell_{q-1}^{\betain} \le \lambda_{q-1}^{-\betain}$, then 
writing $\betain = \theta\beta$ for a $\theta>1$, it suffices to take  $\theta<b^3$, and then $\alpha\ll 1$. Choosing $b=6$, this gives the relation 
\begin{align}
    \beta< \frac{\betain}{216}. \label{e:constraint-vinreg-wdq-theta}
\end{align}  The $H^1$ estimate follows because
$ \|\nabla \wdq \|_{L^2} \le  \|\vin-\vin*\psi_{\ell_{q-1}}\|_{L^2}
    \ell_q^{-1}$.
So we have \eqref{e:wq+1-L2}  by summing the  estimates together and taking $a\gg1$.
 
Lastly, the rather crude $L^4$ estimate \eqref{e:wq+1-L4-crude}  follows similarly.
\end{proof}

\section{The new stress error}
\label{s:new-stress}
By the construction of $v_{q+1}=\vv_q+w_{q+1}$ and the equation \eqref{e:subsol-glued-euler} for $(\vv_q,\pp_q,\RRR_q)$, \begin{align*}
    \del_t v_{q+1} + \div(v_{q+1} \otimes v_{q+1})
     & = \del_t \vv_q + \div (\vv_q \otimes \vv_q)  + \nabla \pp_q  - \div\RRR_q \\
    & \qquad + \div \RRR_q + \div(w_{q+1}\otimes w_{q+1}) - \nabla \pp_q
    \\
    & \qquad + \del_t w_{q+1} + \vv_q \cdot \nabla w_{q+1} + \div(w_{q+1} \otimes \vv_q).
\end{align*}
The first line of the right-hand side vanishes by \eqref{e:subsol-glued-euler}. By defining 
\begin{align}
    \Rosc &\coloneq  
    \mathcal R\left(  \div (w_{q+1}\otimes w_{q+1}+\RRR_q)
    -\nabla \bigg(\sum_{j=0}^{\jmax}\sum_{i=0}^{\imax} \eta^2_{ i}\chi^2_\bj\rho_{q,j}\bigg)\right) \label{e:Rosc},\\
\Rlinear&\coloneq \mathcal R( \del_t w_{q+1}) +( w_{q+1}  \ootimes \vv_q+ \vv_q \ootimes w_{q+1}) \notag
\\
p_{q+1}&\coloneq \pp_q + \sum_{j=0}^{\jmax}\sum_{i=0}^{\imax} \eta^2_{ i}\chi^2_\bj\rho_{q,j} + 2w_{q+1}\cdot \vv_q , \label{e:new-pressure}\\
 \RR_{q+1} &\coloneq  \Rosc +\Rlinear\notag,
\end{align} 
it follows that $(v_{q+1},p_{q+1}, \RR_{q+1})$ solves the relaxed Euler system \eqref{e:subsol-euler} with $q$ replaced by $q+1$. In addition, $\Rosc$ and $\Rlinear$ are both symmetric and trace-free.



\subsection{Estimates of the stress error}
Below, we will use the parameter estimates
\begin{align}
 \delta_{q+1}^{1/2}\lambda_{q+1}^{-1/6} 
 + \tau_q^{-1} \ell_q^{-3} \lambda_{q+1}^{-7/6+\alpha/3}
 + \frac{\ell_q^{-4}\delta_{q+1}}{\lambda_{q+1}^{5/6}}
 \lesssim   \lambda_{q+1}^{-4\alpha}\delta_{q+2}.
 \label{e:param-estimates-for-estimating-Rq+1}
 \end{align}
 Each term on the left is bounded by the right hand-side after taking $\alpha\ll1$ if we impose the following parameter restrictions 
\begin{align}
-\beta-\frac16<-2\beta b 
, \notag
\\
\frac{-\beta + 1 + \frac92(-1-\beta(b-1))}b -\frac76 < -2\beta b
, \notag
\\
\frac{4((-1-\beta(b-1))}b -\frac56 -2\beta < -2\beta b
.\notag
\end{align}
Taking $b=6$, we get
\begin{align}
 \beta <\min\Big(\frac1{66},\frac{21}{106},\frac9{40}\Big) =\frac1{66}. \label{e:constraint-b-beta-strong} \end{align}
\begin{prop}[Estimate for $\Rlinear $]\label{p:Rlinear-est}
\begin{align}\notag
\|\Rlinear \|_{L^1}  \lesssim\delta_{q+1}^{1/2}\lambda_{q+1}^{-1/6}
+\tau_q^{-1} \ell_q^{-3} \lambda_{q+1}^{-7/6+\alpha/3}. \end{align}
\end{prop}
\begin{proof}
Lemma \ref{p:improved-holder}  is applicable with constants $p=1,\mu = \ell_q^{-1}$ and $\kappa = \lambda_q^{5/6}$ under the condition \eqref{e:constraint-b}. Hence, 
we deduce that
\begin{align*}
\|w_{q+1}  \ootimes \vv_q+ \vv_q \ootimes w_{q+1}\|_{L^1}
    &\lesssim   \|\vv_q\|_{L^2}  \sup_{i\ge 0} \sum_{k\in\Lambda_i} \bigg\| \sum_{j} a_\bk \bigg\|_{L^2}\|W_\bk\|_{L^1}\\
&\lesssim \delta_{q+1}^{1/2}\lambda_{q+1}^{-1/6},
\end{align*}
The second term is dealt with easily as $V_\bk$ does not depend on $t$:
\begin{align*}
\MoveEqLeft \|\mathcal R( \del_t w_{q+1})\|_{L^1}
=  \frac{1}{\lambda_{q+1}}\bigg\|\mathcal R\curl\sum_{j,i,k} (\del_t a_\bk) V_\bk \bigg\|_{L^{1+\alpha}}\\
&\lesssim \frac{1}{\lambda_{q+1}}  \sup_{i\ge 0} \sum_{k\in\Lambda_i} \bigg\| \sum_{j} \del_t  a_\bk \bigg \|_{L^{1+\alpha}}\|V_\bk \|_{L^{1+\alpha}} \\
&\lesssim   \tau_q^{-1}\ell_q^{-3} \lambda_{q+1}^{-1} (\lambda_{q+1}^{1/6})^{2(\frac12-\frac1{1+\alpha})} \\
&\le \tau_q^{-1} \ell_q^{-3} \lambda_{q+1}^{-7/6+\alpha/3}.\qedhere 
\end{align*}
\end{proof}

\begin{prop}[Estimate for $\Rosc$]\label{p:Rosc-est}
\begin{align} \|\Rosc\|_\alpha  
\lesssim
    \ell_q^{-4}\lambda_{q+1}^{-5/6}\delta_{q+1}
    +\lambda^{-4\alpha}_{q+1}\delta_{q+2}.  \label{e:Rosc-est}
\end{align}
\end{prop}
\begin{proof}
The principal term of $w_{q+1}\otimes w_{q+1}$ is $\wpq\otimes \wpq$, which cancels $\RRR_q$ up to small error and a gradient: writing $\mathbb P_{>0}f := f - \int_{\mathbb T^3}f(y)\dd y$, 
\begin{align}
&\div( \wpq\otimes\wpq + \RRR_q )\notag\\
&=\div\bigg( \sum_{j,i } \eta^2_{ i}\chi^2_\bj\rho_{q,j} \sum_{k \in \Lambda_i} \gamma_k^2(\TR j) W_\bk\otimes W_\bk    + \RRR_q \bigg)\notag\\
&\overset{\mathclap{\eqref{e:TRj-defn}}}=\  \div\bigg( \sum_{j,i } \eta^2_{ i}\chi^2_\bj\rho_{q,j} \bigg(\sum_{k \in \Lambda_i} \gamma_k^2(\TR j)  W_\bk\otimes W_\bk   - \TR j + \Id\bigg)\bigg)\notag\\
&\overset{\mathclap{\eqref{e:R-decomp}}}=\div\bigg( \sum_{j,i } \eta^2_{ i}\chi^2_\bj\rho_{q,j} \sum_{k \in \Lambda_i} \gamma_k^2(\TR j) \mathbb{P}_{>0}(W_\bk \otimes W_\bk )
+\Id \sum_{j,i } \eta^2_{ i}\chi^2_\bj\rho_{q,j}\bigg) \notag  \\
&=\div\bigg( \sum_{j,i,k } a_\bk^2 \mathbb{P}_{>0}(W_\bk \otimes W_\bk )\bigg )
+\nabla\bigg(\sum_{j,i } \eta^2_{ i}\chi^2_\bj\rho_{q,j}\bigg) \eqcolon \div R'+\nabla P'\notag.
\end{align}
The pressure term $P'$ was specifically transferred from \eqref{e:Rosc} to \eqref{e:new-pressure}. We estimate the remaining terms using the fact \eqref{e:wbk-properties} that $W_\bk$ is a stationary Euler flow, and Lemma \ref{p:improved-oscillatory-decay}, which is applicable with parameters  $\kappa = \lambda_{q+1}^{5/6}, \mu=\ell_q^{-1}$ and $L\gg 1$
for $0<p-1\ll1 $ under the condition \eqref{e:constraint-b}: 
\begin{align}
 \|\mathcal R\div R'\|_{L^1} &\lesssim \left\|\mathcal R  \sum_{j,i,k }\div\Big(  a_\bk^2 \mathbb{P}_{>0}(W_\bk \otimes W_\bk )\Big)\right\|_{L^{p}}\notag\\
&=\left\|\mathcal R \sum_{j,i,k } \Big( \nabla( a_\bk^2) \mathbb{P}_{>0}(W_\bk \otimes W_\bk )\Big)\right\|_{L^{p}}\notag\\
&\overset{\mathllap{\textup{Lem. \ref{p:improved-oscillatory-decay}}}} \lesssim \frac{1}{\lambda_{q+1}^{5/6}}\left \lVert  \sum_{j,i,k}  a_\bk \cdot \nabla a_\bk   \right\|_{L^{p}}\overset{\eqref{e:a-bk-Dn-Lp},\eqref{e:rho-sum-bd}}\lesssim \frac{\ell_q^{-4}\delta_{q+1} }{\lambda_{q+1}^{5/6}}  \notag .
\end{align}
Now, $\Rosc = \mathcal R \div R' + \mathcal R \div( w_{q+1}\otimes w_{q+1} - \wpq\otimes \wpq)$, and for $0<p-1\ll1$ depending only on $b$, $\beta$, and $\alpha$, with  $\theta \ll 1$ such that $\frac1p = \frac{1-\theta}1+\frac\theta2$: \begin{align}
\MoveEqLeft \|\mathcal R\div(w_{q+1}\otimes w_{q+1} - \wpq\otimes \wpq) \|_{L^1}\notag 
\lesssim  \|w_{q+1}\otimes w_{q+1} - \wpq\otimes \wpq \|_{L^p}\notag  \\
&\lesssim   \|w_{q+1}\otimes w_{q+1} - \wpq\otimes \wpq \|_{L^1}^{1-\theta}\|w_{q+1}\otimes w_{q+1} - \wpq\otimes \wpq \|_{L^2}^{1-\theta}\notag  \\
&\lesssim   \big((\|\wcq\|_{L^{2}}+\|\wdq\|_{L^{2}})(\|\wpq\|_{L^{2}}+\|\wcq\|_{L^{2}}+\|\wdq\|_{L^{2}})\big)^{1-\theta} \notag 
\\
&\qquad \times  (\|\wpq\|_{L^{4}}+\|\wcq\|_{L^{4}}+\|\wdq\|_{L^{4}})^{2\theta }  \notag 
\\
&\lesssim  (\ell_q^{-4}\lambda_{q+1}^{-1}\delta_{q+1}+\lambda^{-4\alpha}_{q+1}\delta_{q+2}\delta_{q+1}^{1/2})^{1-\theta} \lambda_{q+1}^{8\theta}\notag%
\\
&\lesssim \ell_q^{-4}\lambda_{q+1}^{-5/6}\delta_{q+1}+\lambda^{-4\alpha}_{q+1}\delta_{q+2}, \label{e:wpow2 - wpqpow2}
\end{align}
This finishes the proof of the claimed $L^1$ estimate \eqref{e:Rosc-est}.
\end{proof}
Using \eqref{e:param-estimates-for-estimating-Rq+1} and $a\gg 1$, Propositions \ref{p:Rlinear-est} and \ref{p:Rosc-est} imply that
 \begin{align*}
    \MoveEqLeft \|\RR_{q+1}\|_{L^1} \le   \|\Rosc\|_{L^1} + \|\Rlinear \|_{L^1}\\
    &\lesssim 
    \ell_q^{-4}\lambda_{q+1}^{-5/6}\delta_{q+1}
    +\lambda^{-4\alpha}_{q+1}\delta_{q+2}
    +\delta_{q+1}^{1/2}\lambda_{q+1}^{-1/6}
    +\tau_q^{-1} \ell_q^{-3} \lambda_{q+1}^{-7/6+\alpha/3}\\
    &\le  \lambda^{-3\alpha}_{q+1}\delta_{q+2},
\end{align*}
which shows that \eqref{e:RR-q-L1} holds with $q+1$ in place of $q$.
\section{Energy iteration}
\label{s:energy-iteration}
In this final section, we show that \eqref{e:energy-q-estimate} holds with $q+1$ in place of $q$, which completes the proof of Proposition \ref{p:main-prop}. This in turn follows from the following proposition.
\begin{prop}[Energy estimate for $v_{q+1}$]
\label{p:energy}For $t\in[1-\tau_q,T]$,
    \begin{align}
        \bigg|e(t)-\int_{\mathbb T^3} |v_{q+1}(t)|^2\dd x -\frac{\delta_{q+2}}2\bigg| \le \frac{\delta_{q+2}}{10}. \label{e:sufficient-energy-estimate}
    \end{align}
\end{prop}
\begin{proof}
Notice that we only want to control the energy for times $t\in[1-\tau_q,T]$. For such times, $\zeta(t)\equiv 0$, so that \eqref{e:rho-q-i} simplifies to
\begin{align}
 \rho_{q,0} (t) 
 &=\frac{ \tilde\rho_{q,0}(t)}{  \sum_{i= i_{q}}^{\imax} \int_{\mathbb T^3}\eta_j^2\chi_{(0)}^2 (\tilde x,t)\dd \tilde x }. \label{e:simpler-rhoq0-for-energy}
\end{align}
In addition, $t<1-\tau_q$ also implies $\eta_{ i}=0$ for $i<\imax$. So all sums in $i$ start from $i_{q}$. We write the total energy gap into the three terms on the right:
\begin{align}\label{e:energy0}
 &e(t)-\int_{\mathbb T^3} |v_{q+1}|^2\dd x-\frac{\delta_{q+2}}2\notag\\
&= \Big(e(t)-\int_{\mathbb T^3} |\vv_q|^2\dd x -\frac{\delta_{q+2}}2\Big) -\int_{\mathbb T^3} |w_{q+1}|^2\dd x- 2 \int_{\mathbb T^3} w_{q+1}\cdot \vv_q \dd x.
\end{align}
The term $\int_{\mathbb T^3} |\wpq|^2 \dd x$ in 
\begin{align*}
 \int_{\mathbb T^3} \! |w_{q+1}|^2\dd x=
 \int_{\mathbb T^3} \! |\wpq|^2+|\wcq|^2+|\wdq |^2
 +2|\wdq \cdot(\wpq+\wcq) |+2|\wpq\cdot \wcq|\dd x
\end{align*}
will cancel $e(t)-\int_{\mathbb T^3} |v_{q+1}|^2\dd x-\frac{\delta_{q+2}}2$ up to small error, due to the definition of $\rho_{q,0}$. All other terms will  be shown to be small.  To this end,
using $\int_{\mathbb T^3} W_\bk\otimes W_\bk \dd x = k\otimes k$ and $\tr \RRR_q=0$, we can rewrite
\begin{align}\label{e:energy2}
    &\int_{\mathbb T^3} |\wpq|^2 \dd x\notag
    = \int_{\mathbb T^3} \tr (\wpq\otimes \wpq) \dd x \notag
    = \int_{\mathbb T^3} \tr\Bigg(\sum_{j,i, k} a_\bk^2 W_\bk\otimes W_\bk \Bigg)  \dd x\notag \\
    &=\int_{\mathbb T^3} \tr\Bigg(\sum_{j,i, k}  a_\bk^2 \mathbb P_{>0}(W_\bk\otimes W_\bk)   \dd x
    +\sum_{j,i} \eta^2_{ i}\chi^2_\bj\sum_{k\in\Lambda_i}\gamma_k^2(\TR j) k\otimes k \Bigg) \dd x\notag \\
    &=\int_{\mathbb T^3} \!\!\!\tr\!\Bigg(\sum_{j,i, k}  a_\bk^2 \mathbb P_{>0}(W_\bk\otimes W_\bk) \Bigg)  \dd x
    +\int_{\mathbb T^3} \!\!\!\tr\!\Bigg(\sum_{j,i, k} \eta^2_{ i}\chi^2_\bj(\rho_{q,j}\Id-\RRR_q) \Bigg)  \dd x\notag \\
     &=\int_{\mathbb T^3} \sum_{j,i, k}  a_\bk^2 \mathbb P_{>0}(|W_\bk|^2)   \dd x +3\int_{\mathbb T^3} \sum_{j,i} \eta^2_{ i}\chi^2_\bj \rho_{q,j}\dd x.
     \notag
\end{align}
By the definition \eqref{e:rho-tilde-q0-defn} of $\tilde \rho_{q,0}(t)$ and \eqref{e:simpler-rhoq0-for-energy}, we deduce that
\begin{align}
3\int_{\mathbb T^3} \sum_{j ,i } \eta^2_{ i}\chi^2_\bj \rho_{q,j} \dd x 
&=3\int_{\mathbb T^3} \sum_{j\ge 1,i} \eta^2_{ i}\chi^2_\bj\dd x 4^{j+j_\gamma}\delta_{q+1}\lambda_q^{-\alpha}+\tilde\rho_{q,0} (t)\notag\\
&=e(t) - \int_{\mathbb T^3} |\vv_q|^2 \dd x - \frac{\delta_{q+2}}2. \notag%
\end{align}
Hence,
\[\int_{\mathbb T^3} |\wpq|^2 \dd x =  \sum_{j,i, k}  \int_{\mathbb T^3}a_\bk^2 \mathbb P_{>0}(|W_\bk|^2)   \dd x +  \left[e(t) - \int_{\mathbb T^3} |\vv_q|^2 \dd x - \frac{\delta_{q+2}}2\right].\]

By `integration by parts' in the form
\[ \left|\int_{\mathbb T^3}f\mathbb{P}_{\geq c}g\dd x\right|
=\left|\int_{\mathbb T^3}|\nabla|^Lf|\nabla|^{-L} \mathbb{P}_{\geq c}g\dd x\right| \le c^{-L}\|f\|_{H^L} \|g\|_{L^2},\]
and the fact that
 $\mathbb{P}_{>0}(|W_\bk|^2)=\mathbb{P}_{\geq {\lambda}_{q+1}^{5/6}}(|W_\bk|^2)$,  we have for all $L\ge0$ that
 \[ \sum_{j,i, k}\left|  \int_{\mathbb T^3}a_\bk^2 \mathbb P_{>0}(|W_\bk|^2)   \dd x \right|\lesssim_L \lambda_{q+1}^{\frac{-5L}6 }\ell_q^{-4L}\delta_{q+1}. \]
 From \eqref{e:constraint-b}, we can assume that $\ell_q^{-4} \lambda_q^{-5/6} \le \lambda_q^{-\alpha_*}$ for a small  positive constant $\alpha_*$ independent of $a$. So, taking $L>\frac{4\alpha b+2\beta b^2}{\alpha_*}$ and $a\gg1$, we obtain
\begin{align}\notag
\ \sum_{j,i, k}\left|  \int_{\mathbb T^3}a_\bk^2 \mathbb P_{>0}(|W_\bk|^2)   \dd x \right|\ll   \lambda_{q+1}^{-4\alpha} \delta_{q+2}. 
\end{align}

Having dealt with the principal part of $w_{q+1}$, it remains to show that the remaining terms in \eqref{e:energy0} are sufficiently small. But
$\int_{\mathbb T^3} 
|\wcq|^2+|\wdq |^2
 +2|\wdq \cdot(\wpq+\wcq) |+2|\wpq\cdot \wcq|
 \dd x \lesssim \lambda_{q+1}^{-4\alpha}\delta_{q+2}$ by \eqref{e:wpow2 - wpqpow2},
 and similarly using the parameter estimates \eqref{e:param-estimates-for-estimating-Rq+1}, and Lemma \ref{p:improved-holder} with $p=1$, $\mu=\ell_q^{-1}$ and $\kappa = \lambda_q^{5/6}$ as in the proof of Proposition \ref{p:Rlinear-est}, we have
 \begin{align}\notag
\left|\int_{\mathbb T^3} w_{q+1}\cdot \vv_q \dd x\right|&\le 
\|\vv_q\|_{L^2}(\|\wdq \|_{L^2}+\|\wcq\|_{L^2})+\Big|\int_{\mathbb T^3} \wpq\cdot \vv_q \dd x\Big|\notag\\
&\lesssim \lambda_{q+1}^{-4\alpha} \delta_{q+2}+\|\vv_q\|_{L^2}\sup_i\sum_{k\in\Lambda_i}\bigg\|\sum_{j} a_\bk  \bigg\|_{L^2}
\|W_\bk \|_{L^1}\notag\\
&\lesssim \lambda_{q+1}^{-4\alpha} \delta_{q+2}+\delta_{q+1}^{1/2}\lambda_{q+1}^{-1/6}
\notag\\
&\lesssim \lambda_{q}^{-3\alpha} \delta_{q+2}. \notag%
\end{align}

Collecting the above estimates together, we obtain 
\[\bigg|e(t)-\int_{\mathbb T^3} |v_{q+1}|^2\dd x -\frac{\delta_{q+2}}2\bigg| \lesssim \lambda_q^{-3\alpha} \delta_{q+2}.\]
The claimed estimate \eqref{e:sufficient-energy-estimate} now follows by choosing $a\gg1$.
\end{proof}

\appendix
\section{Inverse Divergence Operator}
\label{s:inverse-div-op}
We recall  \cite[Def. 1.4]{zbMATH06456007} the following  operator of order $-1$,
\begin{align*}
    \mathcal R u \coloneq -(-\Delta)^{-1} (\nabla u + \nabla u^\TT ) - \frac12(-\Delta)^{-2} \nabla^2 \nabla\cdot u  + \frac12 (-\Delta)^{-1} (\nabla\cdot u )I_{3\times 3},
\end{align*}
where $U=(-\Delta)^{-1}v$ is the mean-free solution to $-\Delta U = u-\dashint_{\mathbb T^3} u.$
 It is a matrix-valued right inverse of the divergence operator for mean-free vector fields, in the sense that
\[ \div\mathcal R u= u - \dashint_{\mathbb T^3} u. \]
In addition, $\mathcal Ru$ is a traceless and symmetric matrix. The operator $\nabla \mathcal R$ is a \CZ operator, so that for $p\in(1,\infty)$,  $
      \| \nabla \mathcal R u\|_{L^p} \lesssim_p \|u\|_{L^p}.  
$
\section{Improved $L^p$ estimates}
We first recall\footnote{Our version differs slightly from the version in \cite{zbMATH07003146} as we were unable to verify precisely the condition stated,  and in any case our torus has a different sidelength, which changes the constants slightly.} \begin{lem}[{\cite[Lemma 3.7]{zbMATH07003146}}] \label{p:improved-holder-BV} Fix integers $K,\kappa,\mu\ge 1$, such that
\[ \frac{\sqrt3 \mu}\kappa\le \frac13, \quad \text{and} \quad \mu^4 \frac{(2\sqrt3\mu)^K}{\kappa ^{K+3}} \le 1.\]
Let $f,g:\mathbb T^3\to\mathbb R$ be smooth functions such that $g$ is $(\mathbb T/\kappa)^3$-periodic, and for a constant $C_f>0$,
\begin{align}
    \|\nabla^N f\|_{L^2} &\le \mymu^N C_f, & N&=1,2,\dots,\myK+4. \notag%
\end{align}
Then,
\[ \|fg\|_{L^2(\mathbb T^3)} \le 5 \|f\|_{L^2(\mathbb T^3)}\|g\|_{L^2(\mathbb T^3)}.\]
\end{lem}
The following lemma is a small modification of \cite[Lemma 2.1]{zbMATH07031840} which itself is a variant of \cite[Lemma 3.7]{zbMATH07003146}, and works for all $p\in[1,\infty]$, but requires larger separation between $\mu$ and $\kappa$ as $p\to\infty$.
\begin{lem}\label{p:improved-holder}
    Fix $p\in[1,\infty]$, and $\mykappa,\mymu  \in\mathbb N$.  There exists a constant $C_p>0$ depending  only on $p$ such that the following holds. Let $f,g:\mathbb T^d \to \mathbb R$ be smooth functions such that $g$ is $(\mathbb T/\mykappa)^d$-periodic%
    , and for a constant $C_f>0$,
\begin{align}
    \|\nabla^N f\|_{L^p} &\le \mymu^N C_f, & N&=1,2,\dots,\left\lceil \frac dp+1\right\rceil+1. \notag
\end{align}
If $\mymu^{d/p+1}< \mykappa^{1/p}$, then up to a constant depending only on $p$, \[ \|fg\|_{L^p(\mathbb T^d)} \lesssim   C_f\|g\|_{L^p(\mathbb T^d)},\] 

\end{lem}
\begin{proof}
This follows directly from \cite[Lemma 2.1]{zbMATH07031840} which shows $ \|fg\|_{L^p} \lesssim_p  \|f\|_{L^p} \|g\|_{L^p} + \kappa^{-1/p} \|f\|_{C^1} \|g\|_{L^p}$, and the  embedding $W^{d/p+1+\epsilon,p}\subset C^1$.    
\end{proof}

Next, we give a variant of \cite[Lemma B.1]{zbMATH07003146}. 
\begin{lem}\label{p:improved-oscillatory-decay}
Let $p\in (1,\infty)$ and $s>3/p$. Also let integers $\mykappa,\mymu,L\in\mathbb N$ be such that
\[ \mykappa \mymu^{s} \left(\frac{2\mymu}\mykappa\right)^L \lesssim 1, \]

 and let $f$ be a $\mathbb{T}^3$-periodic function such that for some constant $C_f$,
\[\|\nabla^jf\|_{L^p}\le  C_f\mymu^j.\]
for all $j\in[1,\myK+4]$. In addition, assume $\int_{\mathbb{T}^3}f(x)\mathbb{P}_{\geq \mykappa}g(x)dx=0$ and $g$ be a $(\mathbb{T}^3/\mykappa)$-periodic function. Then, 
$$\||\nabla|^{-1}(f\mathbb{P}_{\geq \mykappa}g)\|_{L^p}\lesssim  C_f\frac{\|g\|_{L^p}}{\mykappa},$$
 where the implicit constant is universal.
\end{lem}
\begin{proof}
Observe that $ \int_{\mathbb T^3} \mathbb P_{\ge \mykappa/2} f \mathbb P_{\ge \mykappa  } g  = 0 $,   by $\int_{\mathbb{T}^3}f(x)\mathbb{P}_{\geq \mykappa}g(x)dx=0$. Furthermore, $\mathbb P_{\ge \mykappa} g$ is $\kappa$-periodic.  Hence, by the estimates \[ \||\nabla|^{-1}\mathbb{P}_{\geq \mykappa}h\|_{L^p}\lesssim \frac{1}{\mykappa}\|h\|_{L^p}\text{\ \ and \ \ }\||\nabla|^{-1}\mathbb{P}_{\neq0 }h\|_{L^p}\lesssim \|h\|_{L^p},\]
we find that for any $s>3/p$,
\begin{align*}
    \||\nabla|^{-1}(f\mathbb{P}_{\geq \mykappa}g)\|_{L^p}&\leq  \||\nabla|^{-1}((\mathbb P_{< \mykappa/2} f)\mathbb{P}_{\geq \mykappa}g)\|_{L^p}
    +  \||\nabla|^{-1}((\mathbb P_{\ge  \mykappa/2} f)\mathbb{P}_{\geq \mykappa}g)\|_{L^p}
    \\
    &\le \frac{1}{\mykappa/2}\|(\mathbb{P}_{<  \mykappa/2}f)\mathbb{P}_{\geq \mykappa}g\|_{L^p}
    +\|(\mathbb{P}_{\ge \mykappa/2}f)\mathbb{P}_{\geq \mykappa}g\|_{L^p}\\
    &\le \frac{1}{\mykappa/2}(\|f\|_{L^p}
        +\mykappa\|\mathbb{P}_{\geq \mykappa/2}f\|_{L^\infty })\|g\|_{L^p}\\
    &\lesssim_{s,p}\frac{1}{\mykappa}(\|f\|_{L^p}
        +\mykappa\lVert\lvert\nabla\rvert^{s}\mathbb{P}_{\geq \mykappa/2}f\|_{L^p})\|g\|_{L^p}  \\
    &\lesssim_{p,L} \frac{1}{\mykappa}(\|f\|_{L^p}
        +2^L \mykappa^{-L+1}\||\nabla|^{L+s}\mathbb{P}_{\geq \mykappa/2}f\|_{L^p})\|g\|_{L^p}\\
    &\lesssim \frac{1}{\mykappa}\bigg(1
        +\mykappa\mu^s\bigg(\frac{2\mymu}{\mykappa}\bigg)^{L}\bigg)C_f\|g\|_{L^p}\\
    &\lesssim \frac{1}{\mykappa}C_f                                                                                                                                                                                                                                                                                               \|g\|_{L^p}. \qedhere 
\end{align*}
\end{proof}

\section{$L^p$ \emph{a priori} estimates}
For completeness, we give a proof of $L^p$ estimates ($p\in (1,\infty)$) for solutions $\lv,\lw,\lP,\lf  \in C^\infty(\mathbb T^3\times [0,T])$ to 
\begin{align}
    \left\{
    \begin{alignedat}{2}
        \del_t \lv + \lw\cdot\nabla \lv + \nabla \lP &= \lf, \\ 
        \nabla \cdot \lv = 0, \ \nabla \cdot \lw &= 0,\\
        \lv|_{t=0}&=0,
    \end{alignedat}     
    \right. \label{e:abstracted-equation}
\end{align}
satisfying in addition $\int_{\mathbb T^3} \lv \dd x = \int_{\mathbb T^3} \lw \dd x = \int_{\mathbb T^3} \lf  \dd x = 0$, and $\int_{\mathbb T^3} \lP \dd x = 0.$ 
\begin{lem} If $\lv,\lw,\lP$, and $\lf$ are as above, then we have the estimate
\begin{align}
    \|\lv(\cdot,t)\|_{L^p(\mathbb T^3)} \lesssim \int_0^t \| \lv(\cdot,t')\cdot\nabla\lw(\cdot,t') \|_{L^p(\mathbb T^3)} + \|\lf(\cdot,t')\|_{L^p(\mathbb T^3)} \dd t'. 
    \label{e:abstracted-Lp-estimate}
\end{align}    
\end{lem}
\begin{proof}
We test the equation against $|\lv|^{p-2}\lv\in L^\infty(\mathbb T^3\times [0,T])$. By $\div\lv=0$, and the fact that $|\lv|^p$ belongs to the \Holder space $C^{1,\min(1,p-1)}(\mathbb T^3\times [0,T])$, we have the pointwise inequality
\[ \del_t \|\lv\|_{L^p(\mathbb T^3)}^p \le (\|\nabla \lP\|_{L^p(\mathbb T^3)}+\|\lf\|_{L^p(\mathbb T^3)} )\|\lv\|_{L^p}^{p-1}. \]
It follows that from the continuity of $\|\lv\|_{L^p}$ that (one checks subsets of $\lv=0$ and $\lv\neq0$ separately)
\[ \|\lv\|_{L^p(\mathbb T^3)} \le \int_0^t \|\nabla \lP(\cdot,t')\|_{L^p(\mathbb T^3)} + \|\lf(\cdot,t')\|_{L^p(\mathbb T^3)} \dd t'.\]
Taking the divergence of \eqref{e:abstracted-equation} leads to the Poisson equation 
\[ \Delta \nabla \lP = \nabla \div \lf - \nabla \div(\lv\cdot\nabla \lw).\]
We used the identity $\div(\lw\cdot\nabla \lv)=\div(\lv\cdot\nabla \lw)$ which holds for divergence free fields. The inequality \eqref{e:abstracted-Lp-estimate} now follows by the $L^p$-boundedness of \CZ operators.
\end{proof}

\section*{Acknowledgements}
This work was supported by the National Key Research and Development Program of China (No. 2020YFA0712900) and NSFC Grant 11831004.

\printbibliography
\end{document}